\documentstyle[amssymb,10pt]{amsart}
\newtheorem{prop}{Proposition}[section]
\newtheorem{prop:def}{Proposition-Definition}[section]

\newtheorem{lemma}{Lemma}[section]
\newtheorem{thm}{Theorem}[section]

\theoremstyle{remark}
\newtheorem{remark}{Remark}

\begin{document} 

\newcommand{\nc}{\newcommand} \nc{\on}{\operatorname}

\nc{\¦}{{|}}

\nc{\pa}{\partial} 

\nc{\cA}{{\cal A}}\nc{\cB}{{\cal B}}\nc{\cC}{{\cal C}} \nc{\cE}{{\cal E}}
\nc{\cG}{{\cal G}}\nc{\cH}{{\cal H}}\nc{\cX}{{\cal X}}\nc{\cR}{{\cal R}}
\nc{\cL}{{\cal L}} \nc{\cK}{{\cal K}}\nc{\cO}{{\cal O}} 
\nc{\cF}{{\cal F}}\nc{\cM}{{\cal M}} \nc{\cW}{{\cal W}} \nc{\cT}{{\cal T}} 

\nc{\sh}{\on{sh}}\nc{\Id}{\on{Id}}\nc{\Diff}{\on{Diff}}
\nc{\ad}{\on{ad}}\nc{\Der}{\on{Der}}\nc{\End}{\on{End}}
\nc{\res}{\on{res}}\nc{\ddiv}{\on{div}} \nc{\FS}{\on{FS}}
\nc{\card}{\on{card}}\nc{\dimm}{\on{dim}}\nc{\gr}{\on{gr}}
\nc{\Jac}{\on{Jac}}\nc{\Ker}{\on{Ker}} \nc{\Den}{\on{Den}}
\nc{\Imm}{\on{Im}}\nc{\limm}{\on{lim}}\nc{\Ad}{\on{Ad}}
\nc{\ev}{\on{ev}} \nc{\Hol}{\on{Hol}}\nc{\Det}{\on{Det}}
\nc{\Cone}{\on{Cone}} \nc{\pseudo}{{\on{pseudo}}}
\nc{\class}{{\on{class}}}\nc{\rat}{{\on{rat}}}
\nc{\local}{{\on{local}}}\nc{\an}{{\on{an}}}
\nc{\Lift}{{\on{Lift}}}\nc{\Mer}{{\on{Mer}}}\nc{\mer}{{\on{mer}}}
\nc{\lift}{{\on{lift}}}\nc{\diff}{{\on{diff}}}\nc{\Aut}{{\on{Aut}}}
\nc{\DO}{{\on{DO}}}\nc{\Frac}{{\on{Frac}}}\nc{\cl}{{\on{class}}}
\nc{\Fil}{{\on{Fil}}}\nc{\id}{{\on{id}}}\nc{\LBA}{\on{LBA}}
\nc{\Lie}{{\on{Lie}}}\nc{\br}{{\on{br}}}

\nc{\Bun}{\on{Bun}}\nc{\diag}{\on{diag}}\nc{\KZ}{{\on{KZ}}}
\nc{\CB}{{\on{CB}}}\nc{\out}{{\on{out}}}\nc{\Hom}{{\on{Hom}}}
\nc{\FO}{{\on{FO}}}\nc{\Irr}{{\on{Irr}}}\nc{\Vect}{{\on{Vect}}}
\nc{\univ}{{\on{univ}}}\nc{\cop}{{\on{cop}}}\nc{\op}{{\on{op}}}
\nc{\Symm}{{\on{Sym}}}\nc{\Alt}{{\on{Alt}}}

\nc{\al}{\alpha}\nc{\de}{\delta}\nc{\si}{\sigma}\nc{\ve}{\varepsilon}\nc{\z}{\zeta}
\nc{\vp}{\varphi} \nc{\la}{{\lambda}}\nc{\g}{\gamma}\nc{\eps}{\epsilon}
\nc{\La}{\Lambda}\nc{\PsiDO}{\Psi\on{DO}}

\nc{\AAA}{{\mathbb A}}\nc{\CC}{{\mathbb C}}\nc{\KK}{{\mathbb K}}
\nc{\NN}{{\mathbb N}}
\nc{\PP}{{\mathbb P}}\nc{\QQ}{{\mathbb Q}}\nc{\RR}{{\mathbb R}}
\nc{\VV}{{\mathbb V}}\nc{\ZZ}{{\mathbb Z}} 

\nc{\bla}{{\mathbf \lambda}} \nc{\bv}{{\mathbf v}}
\nc{\bz}{{\mathbf z}}\nc{\bt}{{\mathbf t}} 
\nc{\bP}{{\mathbf P}}\nc{\kk}{{\mathbf k}} 
\nc{\Assoc}{{\mathbf Assoc}} 

\nc{\A}{{\mathfrak a}}\nc{\B}{{\mathfrak b}}\nc{\C}{{\mathfrak c}}
\nc{\G}{{\mathfrak g}}
\nc{\HH}{{\mathfrak h}}\nc{\mm}{{\mathfrak m}}\nc{\N}{{\mathfrak n}}
\nc{\T}{{\mathfrak t}}\nc{\xx}{{\mathfrak x}}
\nc{\SG}{{\mathfrak S}}\nc{\grt}{{\mathfrak grt}}

\nc{\unu}{{\underline{\nu}}}\nc{\bl}{{\mathbf l}}

\nc{\wt}{\widetilde}\nc{\ul}{\underline}
\nc{\wh}{\widehat}
\nc{\SL}{{\mathfrak{sl}}}\nc{\ttt}{{\mathfrak{t}}}
\nc{\GL}{{\mathfrak{gl}}}

\nc{\bn}{\begin{equation}}
\nc{\en}{\end{equation}}

\title[A cohomological construction of quantization functors]
{A cohomological construction of quantization functors of Lie bialgebras}

\begin{abstract}
We propose a variant to the Etingof-Kazhdan construction of quantization 
functors. We construct the twistor $J_\Phi$ associated to an associator 
$\Phi$ using cohomological techniques. We then introduce a criterion ensuring 
that the ``left Hopf algebra'' of a quasitriangular QUE algebra is flat. 
We prove that this criterion is satisfied at the universal level. This 
provides a construction of quantization functors, equivalent to the 
Etingof-Kazhdan construction. 
\end{abstract}

\author{B. Enriquez}

\address{IRMA (CNRS), Universit\'e  Louis Pasteur, 7, rue Ren\'e Descartes, 
67084 Strasbourg, France}

\maketitle

\subsection*{Introduction}

In \cite{EK,EK:2}, Etingof and Kazhdan constructed quantization 
functors of Lie bialgebras over a field $\KK$ of characteristic $0$. 
Such a functor is a morphism of 
props $Q : \ul{\on{Bialg}} \to \wh S^\bullet(\wh{\ul{\LBA}})$ with 
suitable classical limit properties. Here $\ul{\on{Bialg}}$ is the 
prop of bialgebras and $\wh{\ul{\LBA}}$ is a completion of the 
prop of Lie bialgebras.  Each quantization functor gives rise 
to a functor $\on{LBA} \to \on{QUE}$ from the category of
Lie bialgebras to that of quantized universal enveloping algebra, inverse to the 
semiclassical limit functor. More precisely, Etingof and Kazhdan 
construct a map $\Phi \to Q_\Phi$ from the set of associators to 
the set of quantization functors.   

The approach to \cite{EK,EK:2} is based on the ``categorical yoga''. 
If $\A$ is a Lie bialgebra, the authors introduce the category of 
$\A$-dimodules, which contains special objects $M_+$ and $M_-^*$. 
They introduce a fiber functor $F : \{\A$-dimodules$\}
\to \KK[[\hbar]]$-modules, with 
$$
F(V) = \Hom_{\A\on{-dimodules}} (M_-,M_+^{*}\otimes V).
$$ 
An important role is played by a morphism $J_\Phi : 
M_-\otimes M_- \to M_-\otimes M_-$, which is a solution to a 
cocycle equation $\wt d(J_\Phi) = \Phi$. 

In this paper, we propose a variant of this construction. The main difference 
of this approach with that of \cite{EK,EK:2} is that (1) we do not use the 
categorical yoga, (2) we construct $J_\Phi$ by cohomological techniques, 
using universal algebras related to the prop $\ul{\on{LBA}}$; moreover, 
we prove that solutions of $\wt d(J_\Phi) = \Phi$ are unique up to gauge. 

Our approach is based on an idea implicit in \cite{EK}, which is formulated 
explicitly in \cite{KT} in the case of finite-dimensional Lie bialgebras. 
Let $\A$ be a finite-dimensional Lie bialgebra, let $\G$ be its double and let 
$r\in \A\otimes\A^* \subset \G^{\otimes 2}$ be the $r$-matrix of $\A$. 
Set $t = r + r^{2,1}$. Then $t\in S^2(\G)^\G$, so 
$(U(\G)[[\hbar]],m_0,\Delta_0,\cR_0 = e^{\hbar t/2}, \Phi(\hbar t^{1,2},\hbar t^{2,3}))$
is a quasitriangular quasi-Hopf algebra. Here $\hbar$ is a formal variable and 
$(m_0,\Delta_0)$ are the product and coproduct of $U(\G)$. Assume that $J\in 
U(\G)^{\otimes 2}[[\hbar]]$ is invertible and is a solution of the 
equation $\wt d(J) = \Phi(\hbar t^{1,2},\hbar t^{2,3})$, where 
$$
\wt d(J) = \big( (1\otimes J)(\id\otimes\Delta_0)(J) \big)^{-1} 
(J\otimes 1)(\Delta_0\otimes \id)(J). 
$$
Then    
$$
(U(\G)[[\hbar]],m_0,\Ad(J) \circ \Delta_0,\cR = J^{2,1}
e^{\hbar t/2} J^{-1}) 
$$
is a quasitriangular Hopf algebra. Then $H_R = \{\langle \cR,\xi\otimes \id \rangle, 
\xi\in U(\G)[[\hbar]]^*\}$ is a Hopf subalgebra of $U(\G)[[\hbar]]$. One 
expects $H_R$ to be a formal series Hopf algebra, whose associated quantized
universal enveloping algebra is a quantization of $\A$. 

The first step of our construction is the solution of the equation 
$\wt d(J_\Phi) = \Phi$, where $J_\Phi$ belongs to a universal 
algebra for solutions of the classical Yang-Baxter equation. 
This step involves cohomology computations (Sections \ref{sect:univ:hoch}
and \ref{sect:map}). We then introduce a criterion ensuring that the 
algebra $H_R$ is flat (Theorem \ref{thm:flatness}). We prove that this 
criterion is satisfied at the universal level (Theorem \ref{thm:form:R}), 
which implies that it is satisfied in the case of each finite-dimensional 
Lie bialgebra.  We derive from there the construction of quantization 
functors. 

In Section \ref{sect:map}, we prove that solutions to the equation 
$\wt d(J_\Phi) = \Phi$ are unique up to gauge. This implies that $J_\Phi$
is gauge-equivalent to Etingof and Kazhdan's twist $J$. Therefore, our
construction is equivalent to theirs.

\medskip 

I would like to express my thanks to P.\ Etingof for numerous discussions 
and collaboration on a project related to this paper.

\section{Universal co-Hochschild cohomology groups} \label{sect:univ:hoch}

\subsection{Schur functors and bifunctors} \label{sect:basic}

Let $n$ be an integer. If $p$ is an idempotent in $\QQ\SG_n$,
we say that $p$ is {\it irreducible} if $p\neq 0$, and $p = p' + p''$, 
with $p',p''$ idempotents, implies $p' = 0$ or $p'' = 0$. We denote by 
$\Irr(n)$ the set of all irreducible idempotents of $\QQ \SG_n$. Two 
idempotents $p,q$ are equivalent, if for some $u\in (\QQ\SG_n)^\times$, 
we have $q = u p u^{-1}$. We write this relation $p \sim q$. 
Then equivalence classes of irreducible idempotents correspond 
bijectively to irreducible finite-dimensional representations of $\SG_n$. 
For each $\bar p\in \Irr(n)/ \sim (= \wh \SG_n)$, we choose a representative $p$ of 
$\bar p$. 

We denote by $\Vect$ the category of vector spaces over $\KK$.  
To a pair $(n,p)$, where $p\in \Irr(n)$, we associate a functor 
$F_{(n,p)} : \Vect \to \Vect$, defined by $F_{(n,p)}(V) = 
p(V^{\otimes n})$. In terms of the representation $S\in \wh \SG_n$
associated to $p$, we have $F_{(n,p)}(V) = (V^{\otimes n} \otimes S)^{\SG_n}$. 
A {\it Schur functor} is determined by a collection of 
vector spaces $(M_{(n,\bar p)})$, where $(n,\bar p)$ runs over all pairs 
of an integer $n$ and of $\bar p\in \Irr(n)/ \sim$; this is the functor defined
by 
$$
F(V) = \bigoplus_{(n,\bar p)} F_{(n,p)}(V) \otimes M_{(n,\bar p)}. 
$$
The spaces $M_{(n,\bar p)}$ will be called the {\it multiplicity spaces}
of the Schur functor $F$.

One defines the direct sum and the tensor product of Schur functors 
by 
$$
(F\oplus G)(V) = F(V) \oplus G(V), \quad 
(F\otimes G)(V) = F(V) \otimes G(V). 
$$

A morphism $\al : F \to G$ between Schur functors 
is the same as a collection of linear maps 
$\al_{(n,\bar p)} : M_{(n,\bar p)} \to N_{(n,\bar p)}$.

The set of all finite Schur functors (i.e., preserving finite-dimensional 
vector spaces) up to isomorphism, identifies with the free abelian semigroup 
$\oplus_{n\geq 0} \NN \wh\SG_n$. The corresponding free group 
$\oplus_{n\geq 0} \ZZ \wh\SG_n$ is then equipped with a commutative ring 
structure. 

In the same way, a {\it Schur bifunctor} is uniquely determined by a collection 
of vector spaces $M_{(n,\bar p),(n',\bar p')}$ (again called the multiplicity 
spaces). The corresponding bifunctor $\Vect^2 \to \Vect$
is defined by 
$$
F(V,W) = \bigoplus_{(n,\bar p),(n',\bar p')}
F_{(n,p)}(V) \otimes F_{(n',p')}(W) \otimes M_{(n,\bar p),(n',\bar p')}. 
$$
One defines the direct sum and tensor product of Schur bifunctors 
by $(F\oplus G)(V,W) = F(V,W) \oplus G(V,W)$ and 
$(F\otimes G)(V,W) = F(V,W) \otimes G(V,W)$. A morphism 
$\al : F \to G$ of Schur bifunctors is a natural 
transformation of bifunctors, and is equivalent to a 
collection of linear maps $\al_{(n,\bar p),(n',\bar p')} : 
M_{(n,\bar p),(n',\bar p')} \to N_{(n,\bar p),(n',\bar p')}$
between multiplicity spaces. 

Moreover, one can define a ``coproduct'' $\Delta$ taking Schur functors to 
bifunctors, by $(\Delta(F))(V,W) = F(V\oplus W)$. $\Delta$ induces a 
bialgebra structure on 
$$
\oplus_{n\geq 0} \ZZ \wh\SG_n. 
$$
One can also define an ``external product'' $\boxtimes$, 
taking a pair $(F,G)$ of Schur functors to the bifunctor $F\boxtimes G$
such that $(F\boxtimes G)(V,W) = F(V)\otimes G(W)$.

\begin{remark}
The Hopf algebra structure on $\oplus_{n\geq 0} \ZZ \wh\SG_n$
may be described as follows. The additive group structure is the 
obvious one. If $S'\in\wh\SG_{n'}$ and $S''\in\wh\SG_{n''}$, we view 
$S'\otimes S''$ as a module over $\SG_{n'}\times \SG_{n''}$, and we set  
$$
S' \cdot S'' = [\on{Ind}_{\SG_{n'} \times \SG_{n''}}^{\SG_{n'+n''}}(S'\otimes S'')]. 
$$
Here the bracket denotes the class in the $K$-group of $\SG_{n'+n''}$, 
which is equal to $\ZZ \wh\SG_{n'+n''}$.  
If $S\in \wh\SG_n$, set 
$$
\Delta(S) = \sum_{n',n''\¦ n' + n'' = n}
[\on{Res}_{\SG_n}^{\SG_{n'} \times \SG_{n''}}(S)], 
$$
where we use the identifications $K_0(\on{Rep}(\SG_{n'} \times \SG_{n''}))
= \ZZ \wh{\SG_{n'} \times \SG_{n''}} = \ZZ \wh\SG_{n'} \otimes 
\ZZ\wh\SG_{n''}$.  

The Littlewood-Richardson coefficients are defined as follows. 
If $n = n'+n''$, $S\in\wh\SG_n$, $S'\in \wh\SG_{n'}$, 
$S''\in \wh\SG_{n''}$, then  
$$
c_{S',S''}^S = [\on{Ind}_{\SG_{n'} \times \SG_{n''}}^{\SG_n}(S' \otimes S'') : S] 
= [\on{Res}_{\SG_n}^{\SG_{n'} \times \SG_{n''}}(S) : S' \otimes S''] 
$$
(the last equality follows from Frobenius reciprocity).  

Then the bialgebra structure is defined by 
$$
S' \cdot S'' = \sum_{S\in\wh\SG_{n'+n''}} c_{S',S''}^S S, 
$$
$$ 
\Delta(S) = \sum_{n',n''\¦ n'+n'' = n} \sum_{S'\in \wh\SG_{n'}, 
S'' \in \wh\SG_{n''}} c_{S',S''}^S S' \otimes S''. 
$$

In \cite{Liu}, it is proved that the graded Hopf algebra 
$\oplus_{n\geq 0} \ZZ \wh\SG_n$ is isomorphic to 
$\ZZ[t_1,t_2,\ldots]$, where each $t_i$ has degree $i$
and with coproduct $\Delta(t_n) = \sum_{i,j\¦ i+j = n} t_i \otimes t_j$. 
The underlying algebraic group is the multiplicative group of invertible
formal series in one variable. 
\end{remark}

\medskip \noindent
{\it Examples.}
We have 
$$
S^n(V\oplus W) = \oplus_{k=0}^n (S^k(V) \otimes S^{n-k}(W)), \; 
\wedge^n(V\oplus W) = \oplus_{k=0}^n (\wedge^k(V) \otimes \wedge^{n-k}(W)) , 
$$
so if $t$ is any scalar, the series $\sum_{n\geq 0} t^n S^n$
and $\sum_{n\geq 0} t^n \wedge^n$ are group-like elements of the completed 
bialgebra $\wh\oplus_{n\geq 0} \ZZ \wh\SG_n$. 

Let us denote by $\otimes^n$ the Schur functor such that $\otimes^n(V) 
= V^{\otimes n}$. Then if $t$ is any scalar, the series $\sum_{n\geq 0}
{t^n\over{n!}} \otimes^n$ is a group-like element of
$\wh\oplus_{n\geq 0} \QQ \wh\SG_n$. 
\medskip

\subsection{Relation to props and operads}

Let $O$ be an operad. So we have a collection of vector 
spaces $O(n)$ with an action of $\SG_n$. Then if $F$ is a 
Schur functor, with $F(V) = \oplus_{(n,\bar p)} F_{(n,p)} \otimes 
M_{(n,\bar p)}$,  we set 
$$
O(F) = \oplus_{(n,\bar p)} p(O(n)) \otimes M_{(n,\bar p)}. 
$$
Then we have $O(F \oplus G) = O(F) \oplus O(G)$, and each morphism 
$\al : F \to G$ of Schur functors induces a linear map $O(\al) : O(F) \to O(G)$. 

Let $P$ be a prop (see \cite{McLane}). So we have a collection of vector spaces 
$P(n,n')$, with an action of $\SG_n \times \SG_{n'}$. Then if $F$
is a Schur bifunctor with multiplicity spaces 
$M_{(n,\bar p),(n',\bar p')}$, with $F(V) = \oplus_{(n,\bar p)} F_{(n,p)} \otimes 
F_{(n',p')} \otimes M_{(n,\bar p),(n',\bar p')}$,  we set 
$$
P(F) = \oplus_{(n,\bar p),(n',\bar p')} (p' \circ P(n,n')\circ p) 
\otimes M_{(n,\bar p),(n',\bar p')}. 
$$
Any morphism $\al : F \to G$ of bifunctors induces a linear map 
$P(\al) : P(F) \to P(G)$. 

If $F$ and $G$ are Schur functors, we set $P(F,G) = P(F \boxtimes G)$. 
Then we have $P(F \oplus F',G) = P(F,G) \oplus P(F',G)$,  
$P(F,G \oplus G') = P(F,G) \oplus P(F,G')$. Any
pair of morphisms $\mu : F \to F'$, $\nu : G \to G'$
induces a bifunctor morphism $\mu \boxtimes \nu : F \boxtimes G 
\to F'\boxtimes G'$ and therefore a linear $P(\mu,\nu) = P(\mu \boxtimes \nu): 
P(F,G) \to P(F',G')$. 

The prop $\ul{\LBA}$ of Lie bialgebras is presented 
by generators and relations mimicking the Lie bialgebra
axioms. Explicitly, the generators are $\mu\in\ul{\LBA}(2,1)$, 
$\delta\in\ul{\LBA}(1,2)$ and the relations are 
$$
\mu\circ(21) + \mu = 0, \; \mu\circ (\mu\boxtimes 1) \circ \big(
(123) + (231) + (312)  \big) = 0,   
$$
$$
(21) \circ\delta + \delta = 0, \;  
\big((123) + (231) + (312)  \big) \circ (\delta\boxtimes 1) \circ \delta = 0,   
$$
$$
\delta \circ \mu = \big((12) - (21)  \big) \circ (1\boxtimes\mu) \circ 
(\delta\boxtimes 1)\circ \big( (12) - (21) \big) . 
$$
We showed (see \cite{Enr:univ,Pos}) that 
$$
\ul{\LBA}(n,m) = \bigoplus_{N\geq 0}
\big( ( (\cF\cL_N)^{\otimes n})_{\sum_i \delta_i} \otimes 
((\cF\cL_N)^{\otimes m})_{\sum_i \delta_i}\big)_{\SG_N} , 
$$
where $\cF\cL_N$ is the free Lie algebra with generators 
$x_i,i = 1,\ldots,N$, graded by $\oplus_{i} \NN \delta_i$
($x_i$ has degree $\delta_i$). 

\subsection{Universal spaces and algebras}

We will construct spaces and algebras, universal for the following 
situation: $\G$ is a Lie bialgebra, $r\in \G \otimes \G$ is a solution of CYBE, 
and of 
$$
(\delta \otimes \id)(r) = [r^{1,3},r^{2,3}], 
\; 
(\id \otimes \delta)(r) = [r^{1,3},r^{1,2}], 
$$
$\G = \A \oplus \B$, $\A,\B$ are Lie subbialgebras of $\G$, 
and $r \in \A \otimes \B$. 

\subsubsection{Universal spaces}
If $F,G$ are Schur functors, we will set 
$$
(F(\A) \otimes G(\B))_\univ = \ul{\LBA}(G,F); 
$$
so we have 
$$
(F(\A) \otimes G(\B))_\univ = 
\bigoplus_{N\geq 1} \big( F(\cF\cL_N)_{\sum_i \delta_i} \otimes 
G(\cF\cL_N)_{\sum_i \delta_i}\big)_{\SG_N} . 
$$

\subsubsection{Universal algebras}

We set 
$$
(U(\G))_\univ = \bigoplus_{N\geq 0} \big( (\cF\cA_N)_{\sum_i \delta_i} 
\otimes (\cF\cA_N)_{\sum_i \delta_i} \big)_{\SG_N}. 
$$
In \cite{quant}, we constructed an algebra structure on 
$(U(\G))_\univ$, such that if $(\G,r = \sum_{i\in I} a_i \otimes b_i)$ 
is a solution of CYBE, the map 
$U(\G)_\univ \to U(\G)$, taking the class of $x_1\cdots x_N \otimes 
y_{\sigma(1)} \cdots y_{\sigma(n)}$ to
$$
\sum_{i_1,\ldots,i_N\in I} a_{i_1} \cdots a_{i_n} b_{i_{\sigma(1)}}
\cdots b_{i_{\sigma(n)}}, 
$$
is an algebra morphism. 

In the same way, one defines 
\begin{equation} \label{def:Ugn}
(U(\G)^{\otimes n})_\univ = \bigoplus_{N\geq 0} \big( 
((\cF\cA_N)^{\otimes n})_{\sum_i \delta_i} 
\otimes ((\cF\cA_N)^{\otimes n})_{\sum_i \delta_i} \big)_{\SG_N} 
\end{equation}
and equips it with an algebra structure, with similar universal properties. 

\subsection{Operations on universal algebras}

\subsubsection{Insertion-coproduct maps}

If $n\leq m$ and $(I_1,\ldots,I_n)$ is a partition of $[1,m]$, 
and if $U(\G)$ is a universal enveloping algebra, define 
$x\mapsto x^{I_1,\ldots,I_n}$ as the linear map $U(\G)^{\otimes n}
\to U(\G)^{\otimes m}$ defined as $x\mapsto ((\Delta^{\¦I_1\¦} \otimes \cdots \otimes 
\Delta^{\¦I_n\¦})(x))^{\sigma}$, where $\sigma$ is the shuffle permutation of 
$[1,m]$, taking $[1,\¦I_1\¦]$ to $I_1$, $\¦I_1\¦+[1,\¦I_2\¦]$ to $I_2$, etc
(here $\¦I\¦$ is the cardinal of $I$ and $[i,j] = \{i,i+1,\ldots,j\}$). 

There are universal maps 
$$
\cop^{I_1,\ldots,I_n} : (U(\G)^{\otimes n})_\univ
\to (U(\G)^{\otimes m})_\univ , 
$$
which are universal versions of the maps $x\mapsto x^{I_1,\ldots,I_n}$. 

\subsubsection{Symmetrization maps}

The symmetrization map $S^\bullet(\G) \to U(\G)$ has no 
normally-ordered version: for example, $r+r^{2,1}$ is 
a normally-ordered element of $U(\G)^{\otimes 2}$, but its image in 
$U(\G)$, $m(r)+m(r^{2,1})$ is not ($m$ is the multiplication map). 

On the other hand, we have a sequence of coalgebra isomorphisms
$$
S^\bullet(\G) \to S^\bullet(\A) \otimes S^\bullet(\B) 
\stackrel{\Symm \otimes \Symm}{\to} U(\A) \otimes U(\B) \stackrel{m}{\to}
U(\G) , 
$$
where the first map is inverse of the composed map 
$S^\bullet(\A) \otimes S^\bullet(\B) \hookrightarrow 
S^\bullet(\G) \otimes S^\bullet(\G) \stackrel{m}{\to} S^\bullet(\G)$. 
This sequence has a universal counterpart, which we now define. 

We construct a linear isomorphism 
\begin{equation} \label{wanted}
(S^\bullet(\A) \otimes S^\bullet(\B))_\univ \to (U(\G))_\univ, 
\end{equation}
as follows: we have 
$$
(S^\bullet(\A) \otimes S^\bullet(\B))_\univ 
= \bigoplus_{N\geq 0} \big( ( S^\bullet(\cF\cL_N) )_{\sum_i \delta_i}
\otimes (S^\bullet(\cF\cL_N) )_{\sum_i \delta_i}\big)_{\SG_N}
$$
Now the symmetrization map of the Lie algebra $\cF\cL_N$
induces a linear isomorphism $S^\bullet(\cF\cL_N) \to U(\cF\cL_N)
= \cF\cA_N$. The tensor square of this map induces the isomorphism 
(\ref{wanted}). 

In the same way, we have natural isomorphisms 
$$
(S^\bullet(\A)^{\otimes n} \otimes S^\bullet(\B)^{\otimes n})_\univ 
\to (U(\G)^{\otimes n})_\univ .
$$
These isomorphisms are universal normally-ordered versions
of the symmetrization isomorphisms. 

\medskip 

{\it Variants.} If $F$ and $G$ are Schur functors, 
we define 
$$
(F(\A) \otimes G(\B) \otimes U(\G))_\univ
$$
as
\begin{equation} \label{space}
\bigoplus_{N\geq 1} \big( 
\big( F(\cF\cL_N) \otimes (\cF\cA_N)^{\otimes n}\big)_{\sum_i \delta_i} 
\otimes 
\big( G(\cF\cL_N) \otimes (\cF\cA_N)^{\otimes n}\big)_{\sum_i \delta_i} 
\big)_{\SG_N} ; 
\end{equation}
if $n\leq m$ and $(I_1,\ldots,I_n)$ is a partition of $[1,m]$, we also have 
linear maps 
$$
(\id_{F(\A) \otimes G(\B)} \otimes \cop^{(I_1,\ldots,I_n)})_\univ :
(F(\A) \otimes G(\B) \otimes U(\G)^{\otimes n})_\univ 
\to (F(\A) \otimes G(\B) \otimes U(\G)^{\otimes m})_\univ, 
$$
which is a universal version of $\id_{F(\A) \otimes G(\B)} \otimes 
\cop^{(I_1,\ldots,I_n)}$. 

A variant of the symmetrization map yields an isomorphism of (\ref{space}) 
with 
$$
\big( F(\A) \otimes S^\bullet(\A)^{\otimes n}
\otimes G(\B) \otimes S^\bullet(\B)^{\otimes n}\big)_\univ. 
$$

If $F$ is a Schur functor, let us write the bifunctor 
$(V,W) \mapsto F(V \oplus W)$ as $\oplus_i F'_i(V) \otimes F''_i(W)$. 
Then we define the space $(F(\G) \otimes U(\G)^{\otimes n})_\univ$
as
$$
\oplus_i \big( F'_i(\A) \otimes F''_i(\B) \otimes U(\G)^{\otimes n}\big)_\univ. 
$$
The map $(\id_{F(\G)} \otimes \cop^{(I_1,\ldots,I_n)})_\univ$ is then 
$\oplus_i (\id_{F'_i(\A) \otimes F''_i(\B)} \otimes \cop^{(I_1,\ldots,I_n)})_\univ$.

\subsection{Universal co-Hochschild complexes}

If $F$ and $G$ are any Schur functors, we define the co-Hochschild 
complex $(C^\bullet_{(F,G)},d^\bullet_{(F,G)})$ as follows. 
We set 
$$
C^n_{(F,G)} = (F(\A) \otimes G(\B) \otimes U(\G)^{\otimes n})_{\univ}, 
$$
and $d^n_{(F,G)} : C^n_{(F,G)} \to C^{n+1}_{(F,G)}$ is equal to 
\begin{align*}
\sum_{k = 1}^n & (-1)^{k+1} (\id_{F(\A) \otimes G(\B)} \otimes 
\cop^{1,2,\ldots,\{k,k+1\},\ldots,n+1})_\univ
\\ & - (\id_{F(\A) \otimes G(\B)} \otimes \cop^{2,\ldots,n+1})_\univ
+(-1)^n (\id_{F(\A) \otimes G(\B)} \otimes \cop^{1,\ldots,n})_\univ. 
\end{align*}

If $F$ is a Schur functor, we set 
$$
(C^\bullet_F,d^\bullet_F) = \oplus_i (C^\bullet_{(F'_i,F''_i)},
d^\bullet_{(F'_i,F''_i)}) ,  
$$
so $C^n_F = (F(\G) \otimes U(\G)^{\otimes n})_\univ$. 
The complexes $(C^\bullet_{(F,G)},d^\bullet_{(F,G)})$, 
resp., $(C^\bullet_{F},d^\bullet_{F})$, are the {\it co-Hochschild 
cohomology complexes} associated with $(F,G)$, resp., with $F$. 

We denote the corresponding cohomology groups by $H^\bullet_{(F,G)}$
and $H^\bullet_F$. 

\subsection{Computation of universal co-Hochschild cohomology}

\begin{thm} \label{thm:coh} \label{thm:1.1}
We have canonical isomorphisms
$$
H^\bullet_{(F,G)} \stackrel{\sim}{\to} (F(\A) \otimes G(\B) 
\otimes \wedge^\bullet(\G))_\univ, 
$$
and 
$$
H^\bullet_{F} \stackrel{\sim}{\to} (F(\G) \otimes \wedge^\bullet(\G))_\univ.  
$$
\end{thm}

To prove this theorem, we first prove several facts about complexes of 
Schur functors and props (Subsections \ref{sect:functors}, \ref{sect:bifunctors}
and \ref{sect:complexes}). Then we will prove Theorem \ref{thm:coh} 
(Subsection \ref{proof:thm:coh})

\subsection{Complexes of Schur functors and operads} \label{sect:functors}

Let us say that a {\it complex of Schur functors} is a  pair 
$(F^\bullet,d^\bullet)$ of (1) a sequence of Schur functors 
$(F^i)_{i = 1,2,\ldots}$, and (2) a collection of morphisms 
$d^i : F^i \to F^{i+1}$ of Schur functors, such that $d^{n+1} \circ d^n = 0$. 

Let us denote by $M^\bullet_{(n,\bar p)}$ the multiplicity space of 
$F^\bullet$  corresponding to $(n,\bar p)$, then each $d^\bullet$
induces a collection of complexes of vector spaces 
$(M^\bullet_{(n,\bar p)},d^\bullet_{(n,\bar p)})$. 

The cohomology of the complex $(F^\bullet,d^\bullet)$ is defined 
as the collection $(H^i(F^\bullet,d^\bullet))_{i= 1,2,\ldots}$ of Schur
functors, where we set 
$$
H^i(F^\bullet,d^\bullet) = \bigoplus_{n\geq 0} \bigoplus_{\bar p\in 
\Irr(n) / \sim} F_{(n,p)} \otimes H^i(M^\bullet_{(n,\bar p)},
d^\bullet_{n,\bar p}). 
$$ 

Let now $O$ be an operad. 
If $(F^\bullet,d^\bullet)$ is a complex of Schur functors, 
then $(O(F^\bullet),O(d^\bullet))$ is a complex of vector spaces. 
Its cohomology can be computed as follows. 

\begin{prop}
For each $i$, there is a canonical isomorphism  
$$
H^i(O(F^\bullet),O(d^\bullet)) \stackrel{\sim}{\to}
O(H^i(F^\bullet,d^\bullet)). 
$$
\end{prop}

{\em Proof.} For each $k$, define $Z^k_{(n,\bar p)}$
as the kernel of $d^k_{(n,\bar p)} : M^k_{(n,\bar p)}
\to M^{k+1}_{(n,\bar p)}$ and $B^k_{n,\bar p}$ as the image of 
$d^{k-1}_{(n,\bar p)} : M^{k-1}_{(n,\bar p)}
\to M^k_{(n,\bar p)}$. Then $B^k_{n,\bar p} \subset Z^k_{(n,\bar p)}$, 
and  $Z^k_{(n,\bar p)} / B^k_{(n,\bar p)} \simeq H^k_{(n,\bar p)}$. 

Let $\wt H^k_{(n,\bar p)}$ be a lift of $H^k_{(n,\bar p)}$ in 
$Z^k_{(n,\bar p)}$; then we have a direct sum decomposition 
$$
Z^k_{(n,\bar p)} = \wt H^k_{(n,\bar p)} \oplus B^k_{(n,\bar p)}.  
$$
Let us also denote by $\Sigma^k_{(n,\bar p)}$ a supplementary of 
$Z^k_{(n,\bar p)}$ in $M^k_{(n,\bar p)}$. Then we have 
$$
M^k_{(n,\bar p)} = \wt H^k_{(n,\bar p)} \oplus B^k_{(n,\bar p)}
\oplus \Sigma^k_{(n,\bar p)}.  
$$
$d^k_{(n,\bar p)}$ vanishes on $\wt H^k_{(n,\bar p)} \oplus B^k_{(n,\bar p)}$, 
and it induces a linear isomorphism $\Sigma^k_{(n,\bar p)} \to 
B^{k+1}_{(n,\bar p)}$. 

We will denote by $K^{k+1}_{(n,\bar p)} : B^{k+1}_{(n,\bar p)}
\to \Sigma^k_{(n,\bar p)}$ the inverse map of $d^k_{(n,\bar p)}$. 
The complex $(M^\bullet_{(n,\bar p)}, d^\bullet_{(n,\bar p)})$
is therefore the direct sum of the complexes $(\wt H^\bullet_{(n,\bar p)},0)$
and $(B^\bullet_{(n,\bar p)} \oplus \Sigma^\bullet_{(n,\bar p)},
(d^\bullet_{(n,\bar p)})_{\¦ B^\bullet_{(n,\bar p)} \oplus \Sigma^\bullet_{(n,\bar p)}})$. 
The first complex has zero differential and therefore
coincides with its cohomology, while the second complex is acyclic 
(a homotopy is $K^\bullet_{(n,\bar p)}$). 

We get 
$$
O(F^\bullet) = \oplus_{(n,\bar p)} p(O(n)) \otimes 
\big( \wt H^\bullet_{(n,\bar p)} \oplus B^\bullet_{(n,\bar p)}
\oplus \Sigma^\bullet_{(n,\bar p)}\big) ;   
$$
this is the sum of the complex
$$
\oplus_{(n,\bar p)} p(O(n)) \otimes 
\wt H^\bullet_{(n,\bar p)} 
$$
with zero differential, and of the complex 
$$
\oplus_{(n,\bar p)} p(O(n)) \otimes 
\big( B^\bullet_{(n,\bar p)}
\oplus \Sigma^\bullet_{(n,\bar p)}\big) ,    
$$
for which $\oplus_{(n,\bar p)} (\id \otimes K^\bullet_{(n,\bar p)})$
is an explicit homotopy, and which is therefore acyclic. \hfill \qed \medskip

\subsection{Complexes of Schur bifunctors and props} \label{sect:bifunctors}

Let us say that a {\it complex of Schur bifunctors} is a pair 
$(F^\bullet,d^\bullet)$ of a collection $(F^i)_{i= 1,2,\ldots}$ of 
Schur bifunctors, and of bifunctor morphisms $d^i : F^i \to F^{i+1}$, 
such that $d^{i+1} \circ d^i = 0$.

Let us denote by $M^\bullet_{(n,\bar p),(n',\bar p')}$ the collection of 
multiplicity spaces of each $F^\bullet$. Then the collection of all
differentials $d^\bullet$ is equivalent to the data of a collection 
$d^\bullet_{(n,\bar p),(n',\bar p')}$ of differentials of each family 
$M^\bullet_{(n,\bar p),(n',\bar p')}$. 

If $P$ is a prop and $(F^\bullet,d^\bullet)$ is a complex of Schur 
bifunctors, we get a complex of vector spaces $(P(F^\bullet),P(d^\bullet))$. 
In the same way as before, we compute the cohomology of this complex: 

\begin{prop} \label{appl:bifunctors}
Let $(F^\bullet,d^\bullet)$ be a complex of Schur bifunctors, 
and let us define its cohomology as the sequence of complexes
$$
H^i(F^\bullet,d^\bullet) = \bigoplus_{(n,\bar p),(n',\bar p')}
F_{(n,p),(n',p')} \otimes H^i(M^\bullet_{(n,\bar p),(n',\bar p')}, 
d^\bullet_{(n,\bar p),(n',\bar p')}). 
$$
Then there is a canonical isomorphism 
$$
H^i(P(F^\bullet),P(d^\bullet))
\stackrel{\sim}{\to}
P(H^i(F^\bullet,d^\bullet)). 
$$
\end{prop}

\subsection{Coproducts of Schur complexes} \label{sect:complexes}

If $F$ is a Schur functor, let us denote by $\Delta(F)$ its 
coproduct; this is a Schur bifunctor, defined by $(\Delta(F))(V,W)
= F(V \oplus W)$ (see Subsection \ref{sect:basic}). $\Delta(F)$ may be expressed in 
terms of multiplicity spaces as follows. Let $M_{(n,\bar p)}$ be the collection of 
multiplicity spaces of $F$. Let us denote by 
$\Delta_{(n',\bar p'),(n'',\bar p'')}^{(n,\bar p)}$ the collection of 
multiplicity spaces of $\Delta(F_{(n,p)})$. Then the multiplicity
space of $\Delta(F)$ corresponding to $((n',\bar p'),(n'',\bar p''))$
is 
$$
\bigoplus_{(n,\bar p)}
\Delta_{(n',\bar p'),(n'',\bar p'')}^{(n,\bar p)}
\otimes M_{(n,\bar p)}. 
$$

Then any morphism $\al : F \to G$ of Schur functors induces a 
morphism $\Delta(\al) : \Delta(F) \to \Delta(G)$ of Schur 
bifunctors: at the level of multiplicity spaces, we have 
$\Delta(\al)_{(n',\bar p'),(n'',\bar p'')} = 
\oplus_{(n,\bar p)} \id \otimes \al_{(n,\bar p)}$. 
We have then the chain rule 
$\Delta(\al \circ \beta) = \Delta(\al) \circ \Delta(\beta)$. 

Then any complex $(F^\bullet,d^\bullet)$ of Schur functors
induces a complex 
$$
(\Delta(F^\bullet),\Delta(d^\bullet))
$$
of Schur bicomplexes. 

\begin{prop} \label{prop:delta}
There is a canonical isomorphism 
$$
H^i(\Delta(F^\bullet),\Delta(d^\bullet))
\to \Delta(H^i(F^\bullet,d^\bullet)). 
$$
\end{prop}

Indeed, both sides identify with 
$$
\bigoplus_{(n',\bar p'),(n'',\bar p'')} 
(F_{(n',p')} \boxtimes F_{(n'',p'')}) \otimes 
\big( \bigoplus_{(n,\bar p)} \Delta_{(n',\bar p'),(n'',\bar p'')}^{(n,\bar p)}
\otimes H^i(M^\bullet_{(n,\bar p)},d^\bullet_{(n,\bar p)})
\big) . 
$$

\subsection{Co-Hochschild complexes at the level of Schur functors}

Let $S^i$ be the ``$i$th symmetric power'' Schur functor and  
$S^\bullet = \oplus_{i\geq 0} S^i$ the ``symmetric algebra'' 
Schur functor. Then we define the co-Hochschild complex of Schur 
functors as follows. It is the sequence 
$$
\cdots \stackrel{d^{n-1}}{\to} 
(S^\bullet)^{\otimes n} \stackrel{d^n}{\to} 
(S^\bullet)^{\otimes n+1} \stackrel{d^{n+1}}{\to} 
\cdots , 
$$
whose specialization to any vector space $V$ is the co-Hochschild 
complex of the free cocommutative coalgebra $S^\bullet(V)$. 

\begin{prop}
The $n$th cohomology Schur functor of this complex is the 
``$n$th alternating power'' Schur functor $\wedge^n$. 
\end{prop}

{\em Proof.} Define $\Psi_n$ as the Schur functor, such that for any
vector space $V$, 
$$
\Psi_n(V) = \sum_{\al = 0}^{n-2} V^{\otimes \al} \otimes S^2(V)
\otimes V^{\otimes (n-2-\al)} \subset V^{\otimes n}. 
$$
Set 
$$
\Psi'_n = \Psi_n \oplus \bigoplus_{(i_1,\ldots,i_n)\neq (1,\ldots,1)}
(S^{i_1} \otimes \cdots \otimes S^{i_n}).  
$$
Then we have $\otimes^n = \wedge^n \oplus \Psi_n$, 
so $(S^\bullet)^{\otimes n} = \wedge^n \oplus \Psi'_n$. 
Now $(d^n)_{\¦ \wedge^n} = 0$ and $d^n$ maps $\Psi'_n$ to $\Psi'_{n+1}$. 
Moreover, $(\Psi'_n,d^n)$ is an acyclic complex (a homotopy can be 
written explicitly). 

{\em Another proof.} Schur functors are faithful on $\Vect$: any  
isomorphism $F \to G$ of Schur functors (i.e., natural transformation of 
functors $\Vect \to \Vect$ with inverse) induces isomorphisms between the 
multiplicity spaces of $F$ and $G$. Since the cohomology of the co-Hochschild
complex $(S^\bullet(V)^{\otimes n}, d^n)$ is equal to $\wedge^n(V)$ for any
$V$, the $n$th cohomology of the complex of Schur functors must be equal to 
$\wedge^n$. 
\hfill \qed \medskip 

\subsection{Coproduct of the co-Hochschild complex}

We have $\Delta(S^\bullet) = S^\bullet \boxtimes S^\bullet$, so the coproduct of the 
co-Hochschild complex has the form 
$$
\cdots  \stackrel{\Delta(d^{n-1})}{\to}
(S^\bullet)^{\otimes n} \boxtimes (S^\bullet)^{\otimes n}
\stackrel{\Delta(d^n)}{\to} (S^\bullet)^{\otimes n+1} \boxtimes 
(S^\bullet)^{\otimes n+1} \cdots . 
$$

We have also $\Delta(\wedge^\bullet) = \wedge^\bullet \boxtimes \wedge^\bullet$, 
so according to Proposition \ref{prop:delta}, we get 
$$
\bigoplus_{n\geq 0} H^n((S^\bullet)^{\otimes n} \boxtimes (S^\bullet)^{\otimes n}, 
\Delta(d^n) ) = \bigoplus_{p,q\geq 0} \wedge^p \boxtimes \wedge^q. 
$$
More generally, if $F$ and $G$ are any Schur functors, we get 
\begin{align} \label{coh:comp}
& \nonumber \bigoplus_{n\geq 0} H^n \big( (F \otimes (S^\bullet)^{\otimes n}) \boxtimes 
(G \otimes (S^\bullet)^{\otimes n}), 
(\id_F \boxtimes \id_G) \otimes \Delta(d_n) \big) 
\\& = \bigoplus_{p,q\geq 0} (F \otimes \wedge^p) \boxtimes (G \otimes \wedge^q). 
\end{align}

\subsection{Proof of Theorem \ref{thm:coh}}
\label{proof:thm:coh}

It suffices to apply Proposition \ref{appl:bifunctors} to $P = \ul{\LBA}$, 
using the cohomology computation (\ref{coh:comp}). \hfill \qed \medskip 

\subsection{The inclusions $\T_n \hookrightarrow T_n \hookrightarrow
(U(\G)^{\otimes n})_\univ$}

In \cite{Dr:gal}, Drinfeld introduced a Lie algebra $\T_n$,   
associated to each integer $n$. Its universal enveloping
algebra is the 
algebra $T_n$, which was introduced in \cite{BN} and called the 
algebra of chord diagrams. $T_n$ is the algebra with generators 
$t_{ij}$, with $i,j$ such that $1\leq i\neq j\leq n$, and relations
$t_{ij} = t_{ji}$, 
\begin{equation} \label{tij:1}
[t_{ij},t_{kl}] = 0
\end{equation}
if $i,j,k,l$ are all distinct, and 
\begin{equation} \label{tij:2}
[t_{ij},t_{ik} + t_{jk}] = 0
\end{equation}
if $i,j,k$ are all distinct. $\T_n$ is the Lie algebra 
with the generators $t_{ij},1\leq i\neq j\leq n$, and 
the same relations. So we have 
$$
T_n = U(\T_n) . 
$$
(The algebras $\T_n$ and $T_n$ are denoted $\A_n^\KK$ and 
$\cA_n^{pb}$ in \cite{Dr:gal} and \cite{BN} respectively.)

The family of algebras $T_n$ is equipped with a system of 
coproduct-insertion operations, which are defined as follows. 
If $n\leq m$ and $I_1,\ldots,I_n$ is a family of disjoint subsets 
of $\{1,\ldots,m\}$, then there is a unique algebra morphism 
$\Delta^{(I_1,\ldots,I_m)}_{T_n\to T_m}$ such that for any pair
$(i,j)$, 
$\Delta^{(I_1,\ldots,I_m)}_{T_n\to T_m}(t_{ij}) = 
\sum_{\al\in I_i,\beta\in I_j} t_{\al\beta}$.

The purpose of the next proposition is to construct a natural 
morphism from $T_n$ to $(U(\G)^{\otimes n})_\univ$, and to prove 
at the same time its injectivity and a basis result for $T_n$. 

For this, we construct elements in $(U(\G)^{\otimes n})_\univ$. 
We define $r$ as the element of $(\A\otimes \B)_\univ$
corresponding to $1_1 \in \ul{\LBA}(1,1) = (\A\otimes \B)_\univ$. 
Then $r$ and $r^{2,1}$ belong to $(\G\otimes \G)_\univ$, 
and we set $t = r + r^{2,1}$. Then for each $i,j$ such that 
$1\leq i\neq j \leq n$, $t^{i,j} \in (U(\G)^{\otimes n})_\univ$.

\begin{prop} \label{piaf}
1) There are unique algebra morphisms $\mu_n : T_n \to (U(\G)^{\otimes n})_{\univ}$, 
such that $\mu_n(t_{ij}) = t^{i,j}$. These morphisms are
compatible with the coproduct-insertion maps, i.e.\ 
$\mu_m \circ \Delta^{(I_1,\ldots,I_m)}_{T_n\to T_m} = 
\Delta^{(I_1,\ldots,I_m)}_{(U(\G)^{\otimes n})_{\univ}
\to (U(\G)^{\otimes m})_{\univ}} \circ \mu_n$. 

2) Each morphism $\mu_n$ is injective. 

3) If $v_1,\ldots,v_n$ are letters, let us denote by 
$F(v_1,\ldots,v_n)$ the free algebra with generators 
$v_1,\ldots,v_n$. Let $(u_{ij})_{1\leq i<j\leq n}$
be letters and let us define a linear map 
$$
\la_n : \bigotimes_{j=2}^n F(u_{ij}, i = 1,\ldots,j-1) \to T_n
$$
by the condition that 
$$
\la_n(\otimes_{j = 2}^n P_j(u_{ij},i = 1,\ldots,j-1))
= P_2(t_{12}) \cdots P_n(t_{1n},\ldots,t_{n-1,n}) . 
$$ 
Then $\la_n$ is a linear isomorphism. 
\end{prop}

{\em Proof.} The proof of 1) is obvious. Let us show that $\la_n$ is
surjective.  By induction, the surjectivity of $\la_n$ may be reduced
to the surjectivity of the linear map $\la_{n-1,n} : T_{n-1} \otimes
F(u_{i1},\ldots,u_{n-1,n}) \to T_n$ sending $P(t_{ij},1\leq i<j < n)
\otimes Q(u_{\al n},\al = 1,\ldots,n-1)$ to $P(t_{ij},1\leq i<j < n)
Q(t_{\al n},\al = 1,\ldots,n-1)$.

Let us prove the surjectivity of $\la_{n-1,n}$. We have the identities
$$
t_{\al n}t_{\beta\gamma} = t_{\beta\gamma} t_{\al n}
$$
if $\al,\beta,\gamma$ are distinct and $<n$, and 
$$
t_{\al n}t_{\al\beta} = t_{\al\beta} t_{\al n} + [t_{\beta n},t_{\al n}]
$$
if $\al,\beta$ are distinct and $<n$. Consider any monomial in the
$t_{ij}$ and apply these identities to a subword of the form $t_{\al
n}t_{\beta\gamma}$, where $\beta,\gamma$ are $<n$. Then in the
resulting monomials, either one $t_{\al n}$ is forwarded to the right,
or the number of letters of the form $t_{\al' n}$ increases
strictly. Applying this procedure repeatedly, we express the initial
monomial as a sum of terms where all $t_{\al n}$ are the the
right. This proves the surjectivity of $\la_{n-1,n}$ and therefore of
$\la_n$.

Let us set 
$$
U^{(n),<} = \bigoplus_{\unu|\nu_{ij} = 0 \on{\ if\ }i\geq j}
(U(\G)^{\otimes n})_{\univ,\unu} 
$$
and let $p : (U(\G)^{\otimes n})_{\univ} \to U^{(n),<}$ be the projection of 
$(U(\G)^{\otimes n})_{\univ}$ to $U^{(n),<}$ parallel to the direct sum of all 
the homogenous components, not appearing the in the decomposition of $U^{(n),<}$. 
Let us set 
$$
F^{(n)} = \bigotimes_{j = 2}^n F(u_{ij},i = 1,\ldots,j-1) , 
$$
and we have 
$F^{(n)} = \oplus_{\kk\in\NN^{\{2,\ldots,n\}}} F^{(n)}_\kk$. 

Introduce the following gradations on $F^{(n)}$ and on $U^{(n),<}$. 
Let $\kk = (k_2,\ldots,k_n)$ be an element of $\NN^{\{2,\ldots,n\}}$. 
If $x_1,\ldots,x_k$ are letters, and if $\ell$ is an integer, let us denote by 
$F(x_1,\ldots,x_k)_\ell$ the subspace of $F(x_1,\ldots,x_k)$ of all elements 
homogeneous of degree $\ell$. Then we set $F^{(n)}_\kk = 
\bigotimes_{j = 2}^n F(u_{ij},i = 1,\ldots,j-1)_{k_j}$. 

On the other hand, for each $\unu$ in $\NN^{\{1,\ldots,n\}^2}$, 
let us denote by $U^{(n),<}_{\unu}$ the homogeneous component of 
$U^{(n),<}$ of degree $\unu$ and let us set 
$$
U^{(n),<}_{\kk} = \bigoplus_{\unu|\forall j, \sum_{i=1}^{j-1} \nu_{ij}
= k_j} U^{(n),<}_{\unu}.
$$
Then we have $U^{(n),<} = \oplus_{\kk\in\NN^{\{2,\ldots,n\}}}
U^{(n),<}_{\kk}$.

The reversed lexicographic order defines a total order on
$\NN^{\{2,\ldots,n\}}$: we say that $\kk \leq \bl$ if either $k_n <
l_n$, or $k_n = l_n$ and $k_{n-1} < l_{n-1}$, etc.

Consider now the sequence of maps
$$
\al_n : F^{(n)} \stackrel{\la_n}{\to} 
T_n \stackrel{\mu_n}{\to} (U(\G)^{\otimes n})_{\univ} \stackrel{p}{\to} 
U^{(n),<}. 
$$
Let us show that  
\begin{equation} \label{inclusion}
\al_n(F^{(n)}_\kk) \subset \bigoplus_{\bl|\bl \leq \kk}
U^{(n),<}_{\bl}.
\end{equation}
Start from a tensor product of monomials $\otimes_{j=2}^n P_j$, where
for each $j$, $P_j\in F(u_{ij},i<j)_{k_j}$.  Replace in this product,
$u_{ij}$ by $t^{i,j}$.  We then replace $t^{i,j}$ by its
value $r^{i,j} + r^{j,i}$ and we put the resulting
expression in the canonical form (in which in each factor, first
components of the $r$-matrix occur before second components). The
maximal expression (for the reversed lexicographic order) in then the
product $\prod_{j=2}^n P_j(r^{i,j},i<j)$; when we reorder all
other expressions using CYBE, the result has strictly smaller degree.

Let us give an example of this transformation. Let us compute the image of 
$1\otimes u_{23}u_{13}$ (it belongs to $F^{(3)}_{(0,0,2)}$). We find 
\begin{align*}
\al_{3}(1\otimes u_{23}u_{13}) =  &  t^{2,3}t^{1,3}
= r^{2,3} r^{1,3} + r^{3,2} r^{1,3}  
+ r^{3,1}r^{2,3} 
\\ & 
+ [r^{2,1},r^{2,3}]+ [r^{3,1},r^{2,1}]
+ r^{3,2} r^{3,1}. 
\end{align*}
The terms of the right side belong to $U^{(3),<}$ and have respectively 
degrees $(0,0,2)$, $(0,1,1)$, $(1,1,0)$, $(1,0,1)$ (first commutator), 
$(2,0,0)$ (second commutator) and $(1,1,0)$. The term of maximal degree 
is the first one. 

If we define a filtration on $F^{(n)}$ and on $U^{(n),<}$ using the
reversed lexicographic order, then the map $\al_n$ is a filtered
map. Its associated graded is the direct sum of the maps $\al_{n,\kk}
: F^{(n)}_{\kk} \to U^{(n),<}_{\kk}$ sending the tensor product
$\otimes_{j = 2}^n P_j(u_{ij},i<j)$ to the product $\prod_{j = 2}^n
P_j(r^{i,j},i<j)$. Since each of these maps is injective, so is
$\al_n$. Therefore $\la_n$ is injective. Since we have already seen
that it is surjective, $\la_n$ is a linear isomorphism. This proves
3).  It then follows that $\mu_n$ is injective, which proves 2).
\hfill \qed \medskip

\begin{remark}
In \cite{Dr:gal}, Section 5, Drinfeld showed the following facts: the
Lie subalgebra of $\T_n$ generated by the $t_{in}, i<n$ is free;
it is an ideal of $\T_n$; $\T_{n-1}$ is a Lie subalgebra of $\T_n$,
so $\T_n$ is the semidirect product of $\T_{n-1}$ by $F_\Lie(\wt
t_{in},i<n)$. So we have a linear isomorphism

\begin{equation} \label{isom:dr}
i_n^{\T} : \T_{n-1} \oplus F_{\Lie}(t_{1n},\ldots,t_{n-1,n})
\to \T_n.
\end{equation}
This isomorphism, together with the isomorphism 
$T_n = U(\T_n)$, implies that the product map 
$T_{n-1} \otimes F(u_{in}, i = 1,\ldots,n-1) \to T_n$
is a linear isomorphism. This gives another proof of 
3) of Proposition \ref{piaf}.  
\hfill \qed\medskip 
\end{remark}

\newpage

\section{A map $\Assoc \to (\wh U_2)^\times$} \label{sect:map}

In what follows, we will denote the algebra $(U(\G)^{\otimes n})_\univ$
by $U_n$. 

\subsection{Algebras $\wh U_n$ and $\wh T_n$}

For each $n\geq 0$, the algebra $U_n = (U(\G)^{\otimes n})_\univ$ is graded: 
its degree $N$ part is 
$$
(U_n)_N = 
\big( (\cF\cA_N^{\otimes n})_{\sum_i \delta_i} \otimes 
(\cF\cA_N^{\otimes n})_{\sum_i \delta_i} \big)_{\SG_N}
$$
(see (\ref{def:Ugn})). 

In the same way, the algebra $T_n$ is graded: we set $\deg(t_{ij}) = 1$
for any $i,j$. 

If $A = \oplus_{N\geq 0} A_N$ is a graded algebra, we denote by 
$\wh A = \wh\oplus_{N\geq 0} A_N$ its completion for the topology defined 
by $(\oplus_{p\geq N} A_p)_{N\geq 0}$. Then the group of invertible elements 
of $\wh A$ is $(\wh A)^\times = \{$elements of $\wh A$, whose component in 
$A_0$ is invertible$\}$. 

This way we introduce algebras $\wh T_n$ and $\wh U_n$. When 
$A = T_n$ or $U_n$, we have $A_0 = \KK$. 

The injection $U_n \hookrightarrow T_n$ is graded, as are the 
insertion-coproduct operations on both sides. So all these operations
extend to the completed algebras. 

\subsection{A noncommutative co-Hochschild complex}

\subsubsection{} For $J \in (\wh{U_2})^\times$, let us define 
$$
\wt d(J) = (J^{2,3} J^{1,23})^{-1} J^{1,2} J^{12,3}. 
$$
So we have a map 
$$
\wt d :  (\wh{U_2})^\times \to (\wh{U_3})^\times. 
$$
 
For $u\in (\wh{U_1})^\times$, $J\in (\wh{U_2})^\times$, 
$\Phi \in (\wh{U_3})^\times$, let us set 
$$
u * J = u^1 u^2 J (u^{12})^{-1}, \; 
u * \Phi = u^{123} \Phi (u^{123})^{-1}
$$
(let us recall that $u^1 = (u\otimes 1)$, $u^{12} = \Delta_0(u)$, etc.). 

This defines actions of $(\wh{U_1})^\times$ on the sets 
$(\wh{U_2})^\times$ and $(\wh{U_3})^\times$. 

\begin{lemma} \label{equiv}
$\wt d$ is $(\wh{U_1})^\times$-equivariant, i.e., 
$$
\wt d(u*J) = u * \wt d(J). 
$$ 
\end{lemma}

\subsubsection{Associators}

Let $\Phi\in (\wh{U_3})^\times$. 
The relation 
\bn \label{duality}
\Phi^{3,2,1} = (\Phi^{1,2,3})^{-1}
\end{equation}
is the duality relation. 
The relation 
\bn \label{pent}
\Phi^{1,2,34} \Phi^{12,3,4} = \Phi^{2,3,4} \Phi^{1,23,4} \Phi^{1,2,3}
\end{equation}
in $(\wh{U_4})^\times$ is called the pentagon relation.  
Finally, the relation 
\bn \label{hex}
e^{t^{2,3} / 2} \Phi^{1,2,3}
e^{t^{1,2} / 2} \Phi^{3,1,2}
e^{t^{1,3} / 2} \Phi^{2,3,1}
= e^{(t^{1,2}+t^{1,3}+t^{2,3})/2}
\end{equation}
in $(\wh{U_3})^\times$ is called the hexagon relation.  

Recall that $(\wh T_3)^\times \hookrightarrow (\wh U_3)^\times$. 
Elements of $(\wh T_3)^\times$ satisfying 
(\ref{duality}), (\ref{pent}), (\ref{hex}) together with 
$$
\Phi_1 = 0, \; \Alt(\Phi_2) = {1\over 8}[t_{12}, t_{23}]
$$
are called 
{\it associators}. Their set is denoted $\Assoc$. 
The set $\Assoc_\Lie$ of {\it Lie associators} is the subset of 
$\Assoc$ of all elements $\Phi$ satisfying the additional relation
$$
\Delta_{\wh{T_3}}(\Phi) = \Phi \otimes \Phi, 
$$
where $\Delta_{\wh{T_3}} : \wh{T_3} \to \wh{T_3} \wh\otimes \wh{T_3}$
is the completion of the universal enveloping algebra coproduct of 
$T_3 = U(\T_3)$. 

(In the definition of $\Assoc$ of \cite{Dr:gal}, (\ref{hex}) is replaced by a pair 
of equations. Each of these equations is equivalent to (\ref{hex}), using 
(\ref{duality}) and the equations $e^{t^{1,2}/2}e^{(t^{1,3}+t^{2,3})/2}
= e^{(t^{1,2}+t^{1,3}+t^{2,3})/2}$ and $e^{t^{2,3}/2}e^{(t^{1,2}+t^{1,3})/2}
= e^{(t^{1,2}+t^{1,3}+t^{2,3})/2}$, which follow from 
the centrality of $t^{1,2}+t^{1,3}+t^{2,3}$. See also \cite{Dr:gal}, Section 5.) 

\subsubsection{A noncommutative co-Hochschild complex}

The infinitesimal version of $\wt d$ and of the map underlying (\ref{pent})
are the co-Hochschild complex differentials.  We have 

\begin{prop} \label{NC:Hoch}
Let $J\in (\wh{U_3})^\times$ and set $\Phi = \wt d(J)$. 
If $\Phi$ is invariant (i.e., belongs to the completion of 
$((U(\G)^{\otimes 3})^\G)_\univ \hookrightarrow U_3$), then 
it satisfies the pentagon equation (\ref{pent}). 
\end{prop}

{\em Proof.} This is a straightforward computation. The key transformation 
if 
\begin{align*}
& \Phi^{2,3,4} \Phi^{1,23,4} \Phi^{1,2,3}
\\ & 
= \Phi^{2,3,4} (J^{1,234})^{-1} (J^{23,4})^{-1} J^{1,23} J^{123,4} \Phi^{1,2,3}
\\ & 
= (J^{1,234})^{-1} \Phi^{2,3,4} (J^{23,4})^{-1} J^{1,23} \Phi^{1,2,3} 
J^{123,4} , 
\end{align*}
where the last equality follows from the fact that $\Phi$ belongs to 
the algebra generated by the $t^{i,j}$ and from Proposition \ref{t:invt}, 
which implies that the $t^{i,j}$ are invariant. 
\hfill \qed \medskip 

\subsection{The main result: the map $\Assoc \to (\wh{U_2})^\times$}

\begin{lemma}
If $J\in (\wh{U_2})^\times$ is such that $\wt d(J)$ is invariant, then for
any $u\in (\wh{U_1})^\times$, we have $\wt d(u*J) = \wt d(J)$. 
\end{lemma}

{\em Proof.} This follows from Lemma \ref{equiv} and Proposition \ref{t:invt}. 
\hfill \qed \medskip 

\begin{thm} \label{thm:main}
There exists a map $\Assoc \to (\wh{U_2})^\times$, $\Phi \mapsto 
J_\Phi$, such that the identity 
$$
\wt d(J_\Phi) = \Phi 
$$
holds for any $\Phi$. 

Moreover, we have 
$$
\{J\in (\wh{U_2})^\times \¦ \wt d(J) = \Phi\} = 
\{u*J_\Phi, u\in (\wh{U_1})^\times\}. 
$$
\end{thm}

Let us add some comments to Theorem \ref{thm:main}. 
If $\al\in\KK$, then 
$$
\wt d(J_\Phi  e^{\al t^{1,2}}) = e^{-\al (t^{1,2} + t^{1,3} + t^{2,3})}
\wt d(J_\Phi)e^{\al (t^{1,2} + t^{1,3} + t^{2,3})} = \Phi,
$$ 
because 
$t^{1,2} + t^{1,3} + t^{2,3}$ is central in $\wh{T_3}$. This is compatible 
with Theorem \ref{thm:main} because we have

\begin{lemma}
Let $m(r)$ be the element of $(\cF\cA_1 \otimes \cF\cA_1)_{\SG_1} 
\subset U(\G)_\univ$ equal to $x_1 \otimes y_1$. 
For any $\al\in \KK$ and $Y\in (U(\G)^{\otimes 2})_\univ$, we have 
$$
Y e^{\al t^{1,2}} = e^{-\al m(r)} * Y. 
$$
\end{lemma}

{\em Proof.} We have $e^{-\al m(r)} * Y = e^{-\al(m(r)^1 + m(r)^2)}
Y e^{\al(m(r)^1 + m(r)^2 + t^{1,2})}$. Now $t^{1,2}$ commutes with $m(r)^1 + m(r)^2$, 
so $e^{-\al m(r)} * Y = e^{-\al(m(r)^1 + m(r)^2)} 
Y e^{\al(m(r)^1 + m(r)^2)} e^{\al t^{1,2}}$. Now Lemma \ref{m(r):central}, 
with $n=2$, says that $[m(r)^1 + m(r)^2,Y] = 0$, which proves the lemma. 
\hfill \qed \medskip 

In \cite{EK,EK:2}, Etingof and Kazhdan constructed a map
$\Phi \mapsto J_\Phi^{\on{EK}}$ from $\Assoc$ to $(\wh{U_2})^\times$, 
such that $\wt d(J_\Phi^{\on{EK}}) = \Phi$. The first part of 
Theorem \ref{thm:main} therefore gives another construction of 
a map with the same properties. It relies on cohomological 
arguments, in contrast with the categorical arguments of 
\cite{EK,EK:2}.  

\subsection{Proof of Theorem \ref{thm:main}}

If $A = \oplus_{i\geq 0} A_i$ is a $\NN$-graded algebra, we denote by 
$A_{\leq n} = \oplus_{i = 0}^n A_i$. We have $\wh A = 
\limm_{\leftarrow} A_{\leq n}$. We denote by $(\wh A^\times)_{\leq n}$
the subset of $A_{\leq n}$ of all elements, whose component in $A_0$
is invertible. So $\wh A^\times = \limm_{\leftarrow} (\wh A^\times)_{\leq n}$. 

Let $\Phi \in \Assoc$. Set $\Phi = 1 + \Phi_1 + \Phi_2 + \cdots$, where 
$\Phi_n \in (T_3)_n$. We will prove inductively the following statement: 

\medskip 
\noindent $(S_n)$ there exists $(J^{(n-1)}_1,\ldots,J^{(n-1)}_{n-1})$, such that 

(a) for each $k$, $J^{(n-1)}_k \in (U_2)_k$,  

(b) $\wt d(1+J^{(n-1)}_1+ \cdots + J^{(n-1)}_{n-1}) = 1 + \Phi_1 + \cdots 
+ \Phi_{n-1}$  (equality in $(U_3)_{\leq n-1}$), 

(c)
\begin{align*}
& \{J\in (U_2)_{\leq n-1} \¦ \wt d(J) = 1 + \Phi_1 + \cdots + \Phi_{n-1}\}
\\ & 
= \{u * (1+J^{(n-1)}_1+ \cdots + J^{(n-1)}_{n-1}) + \la | 
u\in (\wh U_1)^\times, \la\in (\wedge^2(\G))_{\univ,n-1}\}
\end{align*} 
\medskip 

The statements $(S_0)$ and $(S_1)$ are obvious. We find $J^{(0)}_1 = -r/2$. 
Let us assume that $(S_n)$ is true and let us show $(S_{n+1})$. We want to solve the
equation 
\begin{equation} \label{wanted'}
\wt d(1+J_1 + \cdots + J_{n}) = 1 + \Phi_1 + \cdots + \Phi_n
\end{equation}
(equality in $(U_3)_{\leq n}$), where $(J_1,\ldots,J_n)$ are such that for each 
$k$, $J_k\in (U_3)_k$. Assume that $(J_1,\ldots,J_n)$ satisfies (\ref{wanted'}), then  
$(S_n)$(c) implies that there exist $u\in (\wh U_1^\times)_{\leq n-1}$ 
and $\la\in (\wedge^2(\G))_{n-1}$, such that  
$$
1+ J_1 + \cdots + J_{n-1} = u (1+J^{(n-1)}_1 + \cdots + J^{(n-1)}_{n-1}) + \la, 
$$
(equation in $(\wh U_1^\times)_{\leq n}$). 

Let $\wt u$ be a lift of $u$ in $(\wh U_1^\times)_{\leq n}$, set 
$J' = (\wt u)^{-1} * (1+ J_1 + \cdots + J_n)$. Then (\ref{wanted})
is equivalent to   
\begin{equation} \label{new:eq}
\wt d(J') = 1 + \Phi_1 + \cdots + \Phi_n
\end{equation}
(equation in $(U_3)_{\leq n}$), and 
$$
J' = 1 + J^{(n-1)}_1 + \cdots + (J^{(n-1)}_{n-1} + \la) + J'_n, 
$$ 
where $\la\in (\wedge^2(\G))_{\univ,n-1}$ and $J'_n \in (U_2)_n$. 
(\ref{new:eq}) is therefore an equation with unknown $(\la,J'_n)$. 

Our purpose is to show that there exists a pair $(\la^{(n)},J^{(n)}_n)$, 
such that the set of solutions of (\ref{new:eq}) has the form 
$$
\{(\la^{(n)},J^{(n)}_n + d(v) + \la'), \; v\in (U_1)_n, 
\la'\in (\wedge^2(\G))_{\univ,n}\}
$$  
we then set 
$$
J_1^{(n)} = J_1^{(n-1)}, \ldots, J_{n-2}^{(n)} = J_{n-2}^{(n-1)}, 
$$
and 
$$
J_{n-1}^{(n)} = J_{n-1}^{(n-1)} + \la^{(n)}. 
$$
In particular, for each $n$, we have $J^{(n)}_1 = -r/2$. 

Equation (\ref{new:eq}) is equivalent to its degree $n$ part. This degree 
$n$ part has the form 
\begin{equation} \label{key:eq}
d(J'_n) = \Phi_n - \langle J^{(n-1)}_{1}, \ldots, J^{(n-1)}_{n-1} + 
\la\rangle, 
\end{equation}
where for $\al_1 \in (U_2)_1, \ldots, \al_{n-1} \in (U_2)_{n-1}$, 
we denote by $\langle \al_1 \ldots, \al_{n-1} \rangle$ the degree $n$
component of $\wt d(1 + \al_1 + \cdots + \al_{n-1})$. 

We have the identity 
$$
\langle -r/2, \al_2, \ldots,\al_{n-1}+ \la \rangle
= \langle -r/2, \al_2, \ldots,\al_{n-1} \rangle + f(\la), 
$$ 
where 
$$
f : (\wedge^2(\G))_{\univ,n-1} \to (U_3)_n
$$
is the linear map such that 
$$
f(\la) = -{1\over 2} \left( r^{1,2}(\la^{1,3} + \la^{2,3}) 
+ \la^{1,2}(r^{1,3} + r^{2,3}) - r^{2,3}(\la^{1,2} + \la^{1,3})
- \la^{2,3}(r^{1,2} + r^{1,3}) \right) .  
$$

\begin{lemma}
Let us set 
$$
[\![ r,\la ]\!] = [r^{1,2},\la^{1,3}] + [r^{1,2},\la^{2,3}] + 
[r^{1,3},\la^{1,3}] + 
[\la^{1,2},r^{1,3}] + [\la^{1,2},r^{2,3}] + 
[\la^{1,3},r^{1,3}] . 
$$
Then $\la\mapsto [\![r,\la]\!]$ is a linear map $(\wedge^2(\G))_{\univ}
\to (\wedge^3(\G))_{\univ}$.  

Moreover, there exists a 
linear map $f' : (\wedge^2(\G))_{\univ,n-1} \to (U_1)_n$, 
such that we have the identity 
$$
f(\la) = d(f'(\la)) - {1\over 6} [\![r,\la]\!]. 
$$
\end{lemma}

{\em Proof.} One checks that $d(f(\la)) = 0$. Moreover, one computes 
$\Alt(f(\la)) = - {1\over {12}}([\![r,\la]\!] - [\![r,\la]\!]^{3,2,1})$. 
(We include the factor ${1\over {n!}}$ in the definition of $\Alt$.) 
Then Theorem \ref{thm:coh} implies the existence of $\la \mapsto f'(\la)$, 
such that we have 
$$
f(\la) = d(f'(\la)) - {1\over {12}} ( [\![r,\la]\!] - [\![r,\la]\!]^{3,2,1}). 
$$
The Lemma then follows from the fact that $\la\mapsto [\![r,\la]\!]$ maps 
$(\wedge^2(\G))_{\univ}$ to $(\wedge^3(\G))_{\univ}$.  Let us prove this fact. 

If $\la \in (\wedge^2(\G))_{\univ}$, there exists a unique 
$\la'\in (\A \otimes \B)_\univ$, such that $\la = \la' - (\la')^{2,1}$. 
If we set $(\delta \otimes \id)_{\univ}(a) = [r^{1,2},a^{1,3}+a^{2,3}]$
and  $(\id \otimes \delta)_{\univ}(a) = [r^{2,3},a^{1,2}+a^{1,3}]$, then 
$$
[\![r,\la']\!] = (\delta\otimes \id)_\univ(\la') 
- (\id \otimes \delta)_\univ(\la') + [r^{1,3}, (\la')^{2,3} - (\la')^{1,2}]. 
$$
On the other hand, 
$$
[\![r,(\la')^{2,1}]\!] = (\delta\otimes \id)_\univ((\la')^{2,1}) 
- (\id \otimes \delta)_\univ((\la')^{2,1}) + [r^{1,3}, (\la')^{3,2} - (\la')^{2,1}],  
$$
so 
$$
[\![r,\la]\!] = (\delta \otimes \id)_\univ(\la) 
+ (\delta \otimes \id)_\univ(\la)^{2,3,1}
+(\delta \otimes \id)_\univ(\la)^{3,1,2}, 
$$
therefore $[\![r,\la]\!] \in (\wedge^3(\G))_\univ$. 
\hfill \qed \medskip 

Equation (\ref{key:eq}) is then written as
\begin{equation} \label{key:eq'}
d(J'_n) + d(f'(\la)) - {1\over 6} [\![r,\la]\!]
= \Phi_n - \langle -r/2, J^{(n-1)}_2,\cdots, 
J^{(n-1)}_{n-1} \rangle. 
\end{equation}

Let us now solve equation (\ref{key:eq'}). We first prove: 

\begin{lemma} \label{pbm:1}
We have 
$$
d(\Phi_n - \langle -r/2, J^{(n-1)}_2,\cdots, 
J^{(n-1)}_{n-1} \rangle) = 0. 
$$
\end{lemma}

{\em Proof of Lemma.} View $1 + J^{(n-1)}_1 + \cdots + J^{(n-1)}_{n-1}$
as an element of $(U_2)_{\leq n}$ (its degree $n$ part is set to zero), 
which we denote by $J^{(n-1)}$, and set $\Phi' = \wt d(J^{(n-1)})$. 
So $\Phi'\in (U_3)_{\leq n}$. Then 
$$
\Phi' = 1 + \Phi_1 + \cdots + \Phi_{n-1} + \langle J^{(n-1)}_1, 
\cdots , J^{(n-1)}_{n-1}\rangle.  
$$
Now the degree $\leq n-1$ part of $\Phi'$ is invariant, and the degree
zero part of $J^{(n-1)}$ is equal to $1$, so $(\Phi')^{1,2,3}$
commutes with $(J^{(n-1)})^{123,4}$ in $(U_4)_{\leq n}$, therefore
the argument of Proposition \ref{NC:Hoch} applies, and $\wt d(\Phi') = 1$. 

On the other hand, view $1 + \Phi_1 + \cdots + \Phi_{n-1}$ as an element of 
$(U_3)_{\leq n}$ and let us denote by $\langle \Phi_1,\ldots,\Phi_{n-1} \rangle$
the degree $n$ part of $\wt d(1 + \Phi_1 + \cdots + \Phi_{n-1})$; this is 
an element of $(U_4)_n$. Then if $\Theta_n$ is any element of $(U_3)_n$, 
we have 
$$
\wt d(1+\Phi_1 + \cdots + \Phi_{n-1} + \Theta_n) = 
1 + d(\Theta_n) + \langle \Phi_1,\ldots,\Phi_{n-1} \rangle . 
$$
We have $\wt d(\Phi) = 1$. Then if we set $\Theta_n = \Phi_n$, we get 
$$
d(\Phi_n) = - \langle \Phi_1,\ldots,\Phi_{n-1} \rangle . 
$$
On the other hand, since $\wt d(\Phi') = 1$, we have 
$$
d(\langle J^{(n-1)}_1, \cdots , J^{(n-1)}_{n-1}\rangle)
= - \langle \Phi_1,\ldots,\Phi_{n-1} \rangle . 
$$
Comparing both equalities, we get the lemma. 
\hfill \qed \medskip  

Lemma \ref{pbm:1} now implies that there exists a unique $\mu\in 
(\wedge^3(\G))_{\univ,n}$ and $K\in (U_2)_n$, such that 
\begin{equation} \label{anvil}
\Phi_n - \langle J_1^{(n-1)},\ldots,J_{n-1}^{(n-1)} \rangle 
= d(K) + \mu . 
\end{equation}
Equation (\ref{key:eq'}) is then written as 
\bn \label{key:eq''}
d(J'_n) + d(f'(\la)) - {1\over 6} [\![r,\la]\!] = 
d(K) + \mu. 
\end{equation}
 
\begin{lemma}
Equation (\ref{key:eq''}) has a solution $(J'_n,\la)$ iff $\mu$ 
belongs to $\{  - {1\over 6} [\![r,\la]\!] , \la\in 
(\wedge^2(\G))_{\univ,n-1}\}$. 
\end{lemma}

{\em Proof.} If (\ref{key:eq''}) has a solution, then the classes of 
both sides in the co-Hochschild cohomology coincide, therefore 
$\mu = - {1\over 6} [\![r,\la]\!]$. Conversely, if $\mu = 
- {1\over 6} [\![r,\la]\!]$, a solution of (\ref{key:eq''}) is 
$(K - f'(\la),\la)$. \hfill \qed \medskip 

\begin{lemma} \label{key:lemma}
If $\mu$ has the form $- {1\over 6} [\![r,\la]\!]$, then the set of 
all solutions of (\ref{key:eq''}) is the set of all pairs 
$(K - f'(\la) + d(u) + \la', \la)$, where $u\in (U_1)_n$
and $\la' \in (\wedge^2(\G))_{\univ,n}$. 
\end{lemma}

{\em Proof.} The equation $[\![r,\la]\!] = 0$, $\la\in (\wedge^2(\G))_\univ$, 
implies $\la = 0$ (see Proposition \ref{cohom:results}). The result now follows
from 
\begin{align*}
& \{K\in (U_2)_n \¦ d(K) = 0\} 
\\ & = \{ d(u) + \la', u\in (U_1)_n, \la'\in (\wedge^2(\G))_{\univ,n-1}\}. 
\end{align*}
\hfill \qed \medskip 

We now introduce and compute some cohomology groups. 
When $k = 0,\ldots,n-1$, define 
$$
(\id^{\otimes k} \otimes \delta \otimes \id^{\otimes n-k-1})_\univ 
: (\G^{\otimes n})_\univ \to (\G^{\otimes k} \otimes \wedge^2(\G) \otimes 
\G^{\otimes n-k-1})_\univ
$$
by 
$$
a \mapsto [r^{k,k+1},a^{1,\ldots,k,k+2,\ldots,n+1}
+ a^{1,\ldots,k-1,k+1,\ldots,n+1}].  
$$
Then we have a complex $( (\wedge^\bullet(\G))_\univ,\pa^\bullet)$, 
where 
$$
(\wedge^k(\G))_\univ \ni x \stackrel{\pa^k}{\mapsto} 
\Alt((\pa \otimes \id^{\otimes k-1})_\univ(x)) \in 
(\wedge^{k+1}(\G))_\univ . 
$$

We will denote by $H^\bullet( (\wedge^\bullet(\G))_\univ,\pa^\bullet)$ 
the cohomology groups of this complex. These cohomology groups may be 
identified with the Lie coalgebra cohomology groups $H^\bullet(\G,\KK)_\univ$. 

\begin{prop} \label{cohom:results}
1) When $k = 2$ and $\la \in (\wedge^2(\G))_\univ$, we have 
$\pa^2(\la) = [\![r,\la]\!]$.

2) $H^2( (\wedge^\bullet(\G))_\univ,\pa^\bullet) = 0$.  
$H^3( (\wedge^\bullet(\G))_\univ,\pa^\bullet)$ is graded by the degree, 
$H^3( (\wedge^\bullet(\G))_\univ,\pa^\bullet)_N = 0$ when $N\neq 2$, 
and $H^3( (\wedge^\bullet(\G))_\univ,\pa^\bullet)_2$ is 1-dimensional and 
spanned by the class of $[t^{1,2},t^{2,3}]$. 
\end{prop}

In particular, $\Ker(\pa^2) = H^2( (\wedge^\bullet(\G))_\univ,\pa^\bullet) 
= 0$. These results will be proved in Appendix \ref{proof:coh:results}. 
(If $\G$ is a finite-dimensional Lie bialgebra, it is the direct sum of 
$\A$ and $\B$ as a Lie coalgebra, so $H^\bullet(\G,\KK) = H^\bullet(\A,\KK)
\oplus H^\bullet(\B,\KK)$; however, this result cannot be translated
immediately at the ``universal'' level.)

According to Lemma \ref{key:lemma} and to Proposition \ref{cohom:results}, 
in order to prove the induction statement, it remains to prove that 
\bn \label{induction}
\Alt((\delta \otimes \id^{\otimes 2})_\univ(\mu)) = 0. 
\end{equation}

\begin{prop}
We have 
$\Alt((\delta \otimes \id^{\otimes 2})_\univ(\mu)) = 0$. 
\end{prop}

{\em Proof.} Let us set 
$$
J = 1 + J_1^{(n-1)} + \cdots + J^{(n-1)}_{n-1} + K, 
$$
where $K$ is determined by equation (\ref{anvil}). $J$ is an element of 
$(U_2)^\times$. Let us twist the associator $\Phi$ by $J$. 

The resulting associator is 
$$
\Phi' = J^{2,3} J^{1,23} \Phi (J^{1,2} J^{12,3})^{-1}. 
$$
Moreover, equation (\ref{anvil}) implies that $\Phi'$  has the expansion 
\bn \label{exp:Phi'}
\Phi' = 1 - \mu + \sum_{k \geq n+1} \Phi'_k, 
\end{equation}
where we recall that $\mu\in (U_3)_n$, and for each $k$, $\Phi'_k\in (U_3)_k$. 

If $X\in U_k$, let us set 
$$
X^{\wt{1,\ldots,\{\ell,\ell+1\},\ldots,k+1}}
= J^{\ell,\ell+1} X^{1,\ldots,\{\ell,\ell+1\},\ldots,k+1} 
(J^{\ell,\ell+1})^{-1}. 
$$
Then we have the equation 
\bn \label{cocycle:Phi'}
(\Phi')^{\wt{1,2,34}} (\Phi')^{\wt{12,3,4}}
= (\Phi')^{2,3,4} (\Phi')^{\wt{1,23,4}} (\Phi')^{1,2,3}. 
\end{equation}

Moreover, if $X\in (U_k)_n$, we have 
$$
X^{\wt{1,\ldots,\{\ell,\ell+1\},\ldots,k+1}}
= X^{1,\ldots,\{\ell,\ell+1\},\ldots,k+1}
- {1\over 2} [r^{\ell,\ell+1},
X^{1,\ldots,\{\ell,\ell+1\},\ldots,k+1}] +
\sum_{n'>n+1} X_{n'},  
$$
where for each $n'$, $X_{n'}$ has degree $n'$. 

Substitute (\ref{exp:Phi'}) in (\ref{cocycle:Phi'}). 
Consider the homogeneous components of the resulting  equation. 
These equations are trivial in degree $\leq n$. In degree
$n+1$, we obtain
$$
d(\Phi'_{n+1}) + \Alt(-{1\over 2}[r^{1,2},\mu^{1,3,4} + \mu^{2,3,4}]) = 0, 
$$  
in other terms, 
$$
d(\Phi'_{n+1}) -{1\over 2} \Alt((\delta \otimes 
\id^{\otimes 2})_\univ(\mu)) = 0. 
$$
Applying the total antisymmetrization to this identity, we get 
$\Alt((\delta \otimes \id^{\otimes 2})_\univ(\mu)) = 0$. 
\hfill \qed \medskip 

As we already said, the equation 
$$
\Alt((\delta \otimes \id^{\otimes 2})_\univ(\mu)) = 0 , 
$$
together with the facts that $\mu$ has degree $>2$ and
 $H^3((\wedge^\bullet(\G))_\univ,\pa^\bullet)_n = 0$ for $n>2$, 
implies the existence of $\la\in (\wedge^2(\G))_\univ$, such that 
$\mu = -{1\over 6} [\![r,\la]\!]$. We may then apply Lemma 
\ref{key:lemma}, thus proving the induction step. This ends the proof of 
Theorem \ref{thm:main}. 
\hfill \qed \medskip

\newpage 

\section{Construction of quantization functors}

Let $\Phi$ be an associator and let $J$ be such that 
$\wt d(J) = \Phi$. Then $\cR := J^{2,1} \Phi J^{-1}$ is a
solution of the quantum Yang-Baxter equation. We first 
show a result on the form of $\cR$ (Section \ref{form:R}). 
We then show a flatness criterion for Hopf algebras arising from 
solution of the Yang-Baxter equations (Section \ref{sect:flatness}).
We show how to $J$ enables us to construct quantization functors
of Lie bialgebras, first in the finite-dimensional case (Section 
\ref{sect:fd}), then in the general case, where we construct a 
suitable morphism of props (\ref{sect:inf:dim}). 

\subsection{The form of $\cR$} \label{form:R}

\subsubsection{} 

$(U(\G)^{\otimes n})_{\univ}$ is the direct sum of its subspaces 
\begin{equation} \label{subspace}
\big( \bigotimes_{\al=1}^n (\cF\cA_N)_{\sum_{i\in I_\al}\delta_i}
\otimes 
\bigotimes_{\al=1}^n (\cF\cA_N)_{\sum_{i\in J_\al}\delta_i}
\big)_{\SG_N}, 
\end{equation}
where $(I_\al)_{\al = 1,\ldots,n}$ and $(J_\al)_{\al = 1,\ldots,n}$
is a pair of partitions of $\{1,\ldots,N\}$. 

We define a degree on $(U(\G)^{\otimes n})_{\univ}$. This degree 
belongs to $\NN^{2n}$. We denote it by 
$$
(\deg_{(1,a)}, \ldots, 
\deg_{(n,a)}, \deg_{(1,b)}, \ldots, \deg_{(n,b)}).
$$ 
The degree of the subspace (\ref{subspace}) 
is $(\on{card}(I_1),\on{card}(J_1), \ldots, 
\on{card}(I_n),\on{card}(J_n))$. We call $\deg_{(i,a)}$ and 
$\deg_{(i,b)}$ the $a$- and $b$-degrees in the $i$th factor. 
 
This degree behaves with respect to the product as follows. 
Let $x,y$ be homogeneous elements of $(U(\G)^{\otimes n})_{\univ}$. 
Assume that for some $\al$, $\deg_{(\al,b)}(x) = 0$ or 
$\deg_{(\al,a)}(y) = 0$. Then $xy$ is a sum of homogeneous 
components $\sum_i z_i$, such that 
$$
\deg_{(\al,a)}(z_i) \geq \deg_{(\al,a)}(x) + \deg_{(\al,a)}(y)   
$$
and  
$$
\deg_{(\al,b)}(z_i) \geq \deg_{(\al,b)}(x) + \deg_{(\al,b)}(y).    
$$
In other words, if the $\al$th factor of an expression is already 
ordered, ordering the other factors can only increase the $a$-
and $b$- degrees of the $\al$th factor.

\subsubsection{}

In particular, we have 
\begin{align*}
& (U(\G)^{\wh\otimes 2})_{\univ} = 
\wh\bigoplus_{N\geq 0} 
\big( 
\bigoplus_{
(I_2,J_2),(I_2,J_2) \on{\ partitions\ of\ }\{1,\ldots,N\} } 
\\ & 
(\cF\cA_N)_{\sum_{i_1\in I_1}\delta_{i_1}}
\otimes 
(\cF\cA_N)_{\sum_{j_1\in J_1}\delta_{j_1}}
\otimes 
(\cF\cA_N)_{\sum_{i_2\in I_2}\delta_{i_2}}
\otimes 
(\cF\cA_N)_{\sum_{j_2\in J_2}\delta_{j_2}}
\big)_{\SG_N}. 
\end{align*}
$(U(\G)\wh\otimes U(\A))_{\univ}$ (resp., 
$(U(\B)\wh\otimes U(\G))_{\univ}$, etc.) corresponds to the
summands, where $J_2 = \emptyset$ (resp., $I_1 = \null$, etc.)
We define $(U(\G)\otimes U(\A)U(\B)_0)_\univ$ as the completed
sum of all summands with $J_2 \neq\null$. 
Then we have a direct sum decomposition 
$$
(U(\G)^{\wh\otimes 2})_\univ =  
(U(\G)\wh\otimes U(\A))_\univ 
\oplus (U(\G)\wh\otimes U(\A)U(\B)_0)_\univ.  
$$
If $x\in (U(\G)^{\wh\otimes 2})_\univ$, we denote by 
$x_{(U(\G)\wh\otimes U(\A))_\univ }$ and $x_{(U(\G)\wh\otimes U(\A)U(\B)_0)_\univ}$
the components of $x$ corresponding to this decomposition. 

\begin{thm} \label{thm:form:R}
Let $\cR\in (U(\G)^{\wh\otimes 2})_\univ$ be a solution of the 
quantum Yang-Baxter equation (QYBE). Let $\cR_N$ be the homogeneous 
components of $\cR$, and assume that $\cR_0 = 1$, 
$\cR_1 = r \in (\A\otimes\B)_\univ\subset 
(U(\G)^{\wh\otimes 2})_\univ$. Then for each $N\geq 1$, we have 
$$
\cR_N \in (U(\G)\wh\otimes U(\A)U(\B)_0)_\univ 
$$
and $\cR_N\in (U(\A)_0 U(\B)\wh\otimes U(\G))_\univ$.
\end{thm}

{\em Proof.} We will only prove the first statement, since the proof of
the second statement is similar. Let us prove it by induction on $N$. 
It is obvious when $N=1$. Let $N$ be an integer $\geq 2$. Assume 
that we have proved that 
$$
(\cR_1)_{(U(\G)\wh\otimes U(\A))_\univ} = 
\ldots = (\cR_{N-1})_{(U(\G)\wh\otimes U(\A))_\univ} = 0. 
$$
For $x\in (U(\G)^{\wh\otimes 2})_\univ$, we set 
\begin{align*}
[\![r,x]\!] = & [r^{1,2},x^{1,3}] + [r^{1,2},x^{2,3}] + [r^{1,3},x^{2,3}] 
\\ & + [x^{1,2},r^{1,3}] + [x^{1,2},r^{2,3}] + [x^{1,3},r^{2,3}].   
\end{align*}
Then $\cR_N$ satisfies 
\begin{equation} \label{eq:R_N}
[\![r,\cR_N]\!] = 
\sum_{ p,p',p'' | p,p',p''>0, p+p'+p'' = N+1} 
\big( \cR_{p''}^{2,3} \cR_{p'}^{1,3} \cR_p^{1,2} - 
\cR_p^{1,2} \cR_{p'}^{1,3} \cR_{p''}^{2,3} \big).  
\end{equation}

Let us denote by $Y_N$ the r.h.s.\  of (\ref{eq:R_N}). It 
belongs to $(U(\G)^{\wh\otimes 3})_\univ$. As before, 
\begin{align} \label{decomp:3}
& (U(\G)^{\wh\otimes 3})_\univ = \wh\bigoplus_{N\geq 0} \Big( 
\bigoplus_{
(I_1,I_2,I_3),(J_1,J_2,J_3)\on{\ partitions\ of\ }\{1,\ldots,N\}}
\\ & \nonumber 
\big( 
\bigotimes_{\al = 1}^3 (\cF\cA_N)_{\sum_{i_\al\in I_\al} \delta_{i_\al}}
\otimes 
(\cF\cA_N)_{\sum_{j_\al\in J_\al} \delta_{j_\al}}
\big) \Big)_{\SG_N}. 
\end{align}
We denote: 

$\bullet$ by $(U(\G)^{\wh\otimes 2} \wh\otimes U(\A)U(\B)_0)_\univ$ 
the sum of all summands with $J_3 \neq \emptyset$, 

$\bullet$ by $(U(\G)\wh\otimes U(\A)\wh\otimes U(\B)_0)_\univ$ 
the sum of all summands with $J_2 = I_3 = \emptyset$ and 
$J_3\neq\emptyset$, 

$\bullet$ by $(U(\G)\wh\otimes U(\A)\wh\otimes U(\B)_0)^{*,*,1}_\univ$ 
the sum of all summands with $J_2 = I_3 = \emptyset$ and $\card(J_3) = 1$. 

If $x$ is an element of $(U(\G)^{\wh\otimes 3})_\univ$ and 
$S$ is any of these
subspaces of $(U(\G)^{\wh\otimes 3})_\univ$, we denote by $x_S$ the projection of 
$x$ on $S$ parallel to the sum of all other components of 
$(U(\G)^{\wh\otimes 3})_\univ$. 

\begin{lemma} \label{lemma:Y_N}

1) $Y_N\in (U(\G)^{\wh\otimes 2}\wh\otimes U(\A)U(\B)_0)_\univ$, 

2) $(Y_N)_{(U(\G)\wh\otimes U(\A)\wh\otimes U(\B)_0)_\univ^{*,*,1}} = 0$.  

\end{lemma}

{\em Proof of Lemma.} The induction assumption means that for each 
$i = 1,\ldots,N-1$, there exists a unique element $S_i \in 
(U(\A)\otimes U(\G)^{\otimes 2})_\univ$, such that the universal 
version of the linear map 
$$
\big( U(\A)\otimes U(\G)^{\otimes 2}\big) \otimes (\A\otimes \B) 
\to U(\G)^{\otimes 2}, 
$$
$$
(A\otimes B\otimes C) \otimes (x\otimes y) \mapsto Ax B \otimes Cy 
$$
takes $S_i \otimes r$ to $\cR_i$. 

Let $k$ be an integer $\leq N-1$. Then 
$\sum_{q,q'|q,q'>0, q+q' = k} \cR_q^{1,3} \cR_{q'}^{2,3}$ 
is the image of 
$$
\sum_{p,q,q'| p,q,q'>0 \on{\ and\ }q+q' = k} (\cR_q^{1,4}S_{q'}^{2,3,4}) 
\otimes r
$$
by the universal version of the linear map 
$$
\big( U(\G) \otimes U(\A) \otimes U(\G)^{\otimes 2}\big) 
\otimes (\A\otimes \B) \to U(\G)^{\otimes 3}, 
$$
$$
(X\otimes Y\otimes Z\otimes T) \otimes (x\otimes y) 
\mapsto X \otimes Y x Z \otimes Ty. 
$$
This map factors as follows: 
$$
\big( U(\G)\otimes U(\A)\otimes U(\G)^{\otimes 2}\big) 
\otimes (\A\otimes\B) \to U(\G)^{\otimes 3}\otimes \B
\to U(\G)^{\otimes 3}, 
$$
$$
(X\otimes Y \otimes Z\otimes T)\otimes (x\otimes y) 
\mapsto X\otimes YzT \otimes T \otimes y \mapsto
X\otimes YzT \otimes Ty.  
$$
There is a universal version of this factorization. Moreover, 
the image of the universal version 
$$
(U(\G)^{\otimes 3} \otimes \B)_\univ \to (U(\G)^{\otimes 3})_\univ
$$ 
of the last map is contained in $(U(\G)^{\otimes 2}\otimes U(\A)
U(\B)_0)_\univ$. So 
$$
\sum_{q,q'|q,q'\leq N-1, q+q' = k} \cR_q^{1,3} \cR_{q'}^{2,3}
$$ 
belongs to $(U(\G)^{\otimes 2}\otimes U(\A)U(\B)_0)_\univ$. 

Since $(U(\G)^{\otimes 2} \otimes U(\A)U(\B)_0)$ is preserved by 
universal version of the operation of exchanging the two first 
tensor factors,  we also have 
$$
\sum_{p,q| p,q>0, p+q = k} \cR_p^{2,3}\cR_q^{1,3}
\in (U(\G)^{\otimes 2} \otimes U(\A)U(\B)_0)_\univ. 
$$
Now $(U(\G)^{\otimes 2} \otimes U(\A)U(\B)_0)_\univ$ is stable under left and 
right multiplications by elements of the form $x\otimes 1$, where 
$x\in (U(\G)^{\otimes 2})_\univ$. Therefore 
$$
\sum_{\stackrel{p,p',p''|p,p',p''>0}
{p+p'+p'' = N+1}} \cR_p^{1,2} \cR_{p'}^{1,3} \cR_{p''}^{2,3} 
\; \on{and} \;    
\sum_{\stackrel{p,p',p''|p,p',p''>0}{p+p'+p'' = N+1}}  \cR_{p''}^{2,3} 
\cR_{p'}^{1,3} \cR_p^{1,2}
$$
both belong to $(U(\G)^{\otimes 2} \otimes U(\A)U(\B)_0)_\univ$.
This proves 1). 

Let us prove 2). Let us denote by $(U(\G)^{\otimes 3})^k_\univ$ the sum of all 
summands of (\ref{decomp:3}), such that $\card(I_3) + \card(J_3) = k$.
(We call $\card(I_3) + \card(J_3)$ the ``total degree in the third
factor''. Recall that $\card(I_3)$ (resp., $\card(J_3)$) is the 
``$a$-degree (resp., $b$-degree) in the third factor''.)
We will prove that 
$$
\sum_{\stackrel{p,p',p''|p,p',p''>0}
{p+p'+p'' = N+1}} \cR_p^{1,2} \cR_{p'}^{1,3} \cR_{p''}^{2,3} \; \on{and} \;   
\sum_{\stackrel{p,p',p''|p,p',p''>0}{p+p'+p'' = N+1}}  \cR_{p''}^{2,3} 
\cR_{p'}^{1,3} \cR_p^{1,2}
$$
both belong to 
$\oplus_{k\geq 2} (U(\G)^{\otimes 3})^k_\univ$.
The second factor of $\cR_{p'}^{1,3}$ ends with a $b$, which 
will remain in the third factor of $\cR_p^{1,2}\cR_{p'}^{1,3}$. 
So $\cR_p^{1,2}\cR_{p'}^{1,3}\in 
\oplus_{k\geq 2} (U(\G)^{\otimes 3})^k_\univ$. 

Now $\cR_{p''}^{2,3}$ also contains a $b$ at the end of its third factor. 
Consider the product $\cR_p^{1,2}\cR_{p'}^{1,3}\cR_{p''}^{4,3}$. It still contains 
a $b$ at the end of its first factor. Moreover, when we reorder 
the remaining terms of the third factor, the degree of these terms 
cannot drop to zero, so the total degree in the third factor  
is $\geq 2$. After we multiply the second and the 4th factor,
we obtain $\cR_p^{1,2}\cR_{p'}^{1,3}\cR_{p''}^{2,3}$, where we should reorder 
the second factor. However, this operation can only {\it increase}
the total degree in the factors which are already ordered, so 
$\cR_p^{1,2}\cR_{p'}^{1,3}\cR_{p''}^{2,3} \in 
\oplus_{k\geq 2} (U(\G)^{\otimes 3})^k_\univ$. 

In the same way, one proves that  
$\cR_{p''}^{2,3} \cR_{p'}^{1,3}\cR_p^{1,2}
\in \oplus_{k\geq 2} (U(\G)^{\otimes 3})^k_\univ$. 
\hfill \qed\medskip

{\it Sequel of the proof of Theorem \ref{thm:form:R}.}
It follows from Lemma \ref{lemma:Y_N} that 
$$
([\![r,\cR_N]\!])_{(U(\G)\otimes U(\A)\otimes 
U(\B)_0)_\univ^{*,*,1}} = 0. 
$$
Let us now decompose $\cR_N = \cR'_N + \cR''_N$, with 
$$
\cR'_N = (U(\G) \otimes U(\A))_\univ,  
\; \cR''_N = (U(\G) \otimes U(\A)U(\B)_0)_\univ.   
$$ 
We want to prove that $\cR'_N = 0$.  

\begin{lemma} 
We have 
$$
([\![r,\cR''_N]\!])_{(U(\G)\otimes U(\A)\otimes 
U(\B)_0)_\univ^{*,*,1}} = 0. 
$$
\end{lemma}

{\em Proof of Lemma.} We have $[\cR_N^{\prime\prime 1,2},r^{1,3}] 
= \cR_N^{\prime\prime 1,2} r^{1,3} 
- r^{1,3} \cR_N^{\prime\prime 1,2}$. Here only the second term 
is not well-ordered, and this happens only in its first factor. 
The first term  is well-ordered, and since $(\cR''_N)^{1,2}$ 
is a sum of terms with positive $b$-degree in the second factor, 
the same is true for $r^{1,3}(\cR_N'')^{1,3}$. So 
$$
(r^{1,3} \cR_N^{\prime\prime 1,2}
)_{(U(\G)^{\otimes 2} \otimes U(\B)_0)^{*,*,1}_\univ}
= 0 . 
$$
In $\cR_N^{\prime\prime 1,2} r^{1,3}$, only the first factor is ill-ordered. 
After we order this expression, the $b$-degree in the second factor can only 
increase, so it remains positive. So 
$(\cR_N^{\prime\prime 1,2} r^{1,3} 
)_{(U(\G)\otimes U(\A) \otimes U(\B)_0)^{*,*,1}_\univ}
= 0$.
Therefore 
\begin{equation} \label{gath:1}
([\cR_N^{\prime\prime 1,2}, r^{1,3}] 
)_{(U(\G)\otimes U(\A) \otimes U(\B)_0)^{*,*,1}_\univ}
= 0 . 
\end{equation}

We have $[\cR_N^{\prime\prime 1,2},r^{2,3}] 
= \cR_N^{\prime\prime 1,2} r^{2,3} 
- r^{2,3} \cR_N^{\prime\prime 1,2}$. Since 
$r^{2,3} \cR_N^{\prime\prime 1,2}$ is well-ordered, and the 
$b$-degree of the second factor of $\cR_N^{\prime\prime 1,2}$ is 
positive, the $b$-degree of the second factor of $r^{2,3}\cR_N^{\prime\prime 1,2}$ 
is also positive. So $(r^{2,3}\cR_N^{\prime\prime 1,2}
)_{(U(\G)\otimes U(\A) \otimes U(\B)_0)^{*,*,1}_\univ} = 0$.   

Consider now $\cR_N^{\prime\prime 1,2}r^{2,3}$. Write 
$$
\cR_N^{\prime\prime 1,2} = \sum_{\stackrel{\sigma\in\SG_N}{\stackrel{p\in 
\{0,\ldots,N\}}{q\in\{1,\ldots,N\}}}} P_{p,q,\sigma}
a_{\sigma(1)} \cdots a_{\sigma(N-p)} b_1\cdots b_{N-q}
\otimes a_{\sigma(N-p+1)} \cdots a_{\sigma(N)} b_{N-q+1} \cdots b_N. 
$$ 
All the terms in the expansion of $\cR_N^{\prime\prime 1,2} r^{2,3}$
have zero projection on $(U(\G)\otimes U(\A) \otimes U(\G))_\univ$, except 
perhaps 
$$
\sum_{\stackrel{\sigma\in\SG_N}{\stackrel{p\in 
\{0,\ldots,N\}}{q\in\{1,\ldots,N\}}}} P_{p,q,\sigma}
a_{\sigma(1)} \cdots a_{\sigma(N-p)} b_1\cdots b_{N-q}
\otimes a_{\sigma(N-p+1)} \cdots a_{\sigma(N)} b 
\otimes [[b,b_{N-q+1}], \cdots, b_N]. 
$$
Then the $b$-degree in the third factor is $\geq 2$, so 
$(\cR_N^{\prime\prime 1,2}r^{2,3}
)_{(U(\G)\otimes U(\A)\otimes U(\B)_0)_\univ^{*,*,1}} = 0$. Finally, 
\begin{equation} \label{gath:2}
([\cR_N^{\prime\prime 1,2},r^{2,3}])_{(U(\G)\otimes U(\A)\otimes U(\B)_0
)_\univ^{*,*,1}} = 0. 
\end{equation}

Now $[r^{1,2},\cR_N^{\prime\prime 1,3}]\in (U(\G)\otimes \B\otimes U(\G))_\univ$, 
so 
\begin{equation} \label{gath:2,5}
([r^{1,2},\cR_N^{\prime\prime 1,3}])_{(U(\G)\otimes U(\A)\otimes U(\B)_0
)_\univ^{*,*,1}} = 0. 
\end{equation}

Consider now $[\cR_N^{\prime\prime 1,3},r^{2,3}] 
= \cR_N^{\prime\prime 1,3} r^{2,3} - r^{2,3} \cR_N^{\prime\prime 1,3}$.
Since $\cR_N^{\prime\prime 1,3} r^{2,3}$ is well-ordered, the $b$-degree the $b$-degree
of its third factor is the sum of that of $\cR_N^{\prime\prime 1,3}$ (which is 
$>0$) and of $r^{2,3}$ (which is $=1$), so it is $\geq 2$. So  
$(\cR_N^{\prime\prime 1,3}r^{2,3}
)_{(U(\G)\otimes U(\G)\otimes U(\B)_0)_\univ^{*,*,1}} = 0$. 

On the other hand, all the terms in $r^{2,3} \cR_N^{\prime\prime 1,3}$
have zero projection on $(U(\G)^{\otimes 2}\otimes U(\B))_\univ$, except perhaps
$$
\sum_{\stackrel{\sigma\in\SG_N}{\stackrel{p\in 
\{0,\ldots,N\}}{q\in\{1,\ldots,N\}}}} P_{p,q,\sigma}
a_{\sigma(1)} \cdots a_{\sigma(N-p)} b_1\cdots b_{N-q}
\otimes [a_{\sigma(N-p+1)}, \cdots [a_{\sigma(N)},a]] \otimes  
b b_{N-q+1} \cdots b_N,  
$$ 
whose $b$-degree in the third factor is $\geq 2$. So 
$(r^{2,3} \cR_N^{\prime\prime 1,3}
)_{(U(\G)^{\otimes 2}\otimes U(\B)_0)^{*,*,1}_\univ} = 0$. 
So 
\begin{equation} \label{to:use}
([\cR_N^{\prime\prime 1,3}, r^{2,3}]
)_{(U(\G)^{\otimes 2}\otimes U(\B)_0)^{*,*,1}_\univ} = 0,  
\end{equation}
in particular,  
\begin{equation} \label{gath:3}
([\cR_N^{\prime\prime 1,3}, r^{2,3}]
)_{(U(\G)\otimes U(\A)\otimes U(\B)_0)^{*,*,1}_\univ} = 0. 
\end{equation}

(\ref{to:use}) implies  $([r^{1,3},\cR_N^{\prime\prime 2,3}]
)_{(U(\G)^{\otimes 2}\otimes U(\B)_0)^{*,*,1}_\univ} = 0$,
which implies  
\begin{equation} \label{gath:4}
([r^{1,3},\cR_N^{\prime\prime 2,3}]
)_{(U(\G)\otimes U(\A)\otimes U(\B)_0)^{*,*,1}_\univ} = 0. 
\end{equation}

Consider now $[\cR_N^{\prime\prime 2,3}, r^{1,2}] = 
\cR_N^{\prime\prime 2,3}r^{1,2} - r^{1,2} \cR_N^{\prime\prime 2,3}$. 
$\cR_N^{\prime\prime 2,3}r^{1,2}$ is well-ordered and the $b$-degree of its second 
factor is $>0$, so its projection on $(U(\G)\otimes U(\A)\otimes U(\B)_0)^{*,*,1}_\univ$
is zero. Now $r^{1,2} \cR_N^{\prime\prime 2,3}$ is ill-ordered in its second factor only,
so after reordering it, the $a$- and $b$-degrees of its third factor can only increase. 
So the projection of $r^{1,2} \cR_N^{\prime\prime 2,3}$ on 
$(U(\G)\otimes U(\A)\otimes U(\B)_0)^{*,*,1}_\univ$ is zero unless perhaps if 
$(a-$degree in third factor, $b$-degree in third factor$)(\cR_N^{\prime\prime 2,3})
= (0,1)$, i.e., if $\cR_N^{\prime\prime} = \la r$ for some scalar $\la$.
But then $[\la r^{2,3},r^{1,2}]   = \la [r^{1,2},r^{1,3}] + \la[r^{1,3},r^{2,3}]$  
has zero projection on $(U(\G)\otimes U(\A)\otimes U(\B)_0)^{*,*,1}_\univ$. So we have always 
\begin{equation} \label{gath:5}
([r^{1,2},\cR_N^{\prime\prime 2,3}]
)_{(U(\G)\otimes U(\A)\otimes U(\B)_0)^{*,*,1}_\univ} = 0. 
\end{equation}
Summing up (\ref{gath:1})-(\ref{gath:5}), we get the lemma. 
\hfill \qed\medskip 

{\it Sequel of the proof of Theorem \ref{thm:form:R}.}
We have therefore 
$([\![r,\cR'_N]\!])_{(U(\G)\otimes U(\A)\otimes 
U(\B)_0)_\univ^{*,*,1}} = 0$. We now show: 

\begin{lemma} 
$\cR'_N\in (U(\B) \otimes \A)_\univ$. 
\end{lemma}
 
{\em Proof of Lemma.} Let us write 
$$
\cR'_N = \sum_{p\in\{0,\ldots,N\},\sigma\in\SG_N}
P_{p,\sigma} a_{\sigma(1)} \cdots a_{\sigma(p)} b_1\cdots b_N
\otimes a_{\sigma(p+1)} \cdots a_{\sigma(N)}. 
$$
The contributions to $([\![r,\cR'_N]\!])_{
(U(\G)\otimes U(\A)\otimes U(\B)_0)^{*,*,1}_\univ}$ are:  

$\bullet$ from $[\cR^{\prime 1,2}_N,r^{1,3}]$, we get
\begin{equation} \label{contr:1}
\sum_{\stackrel{p\in\{0,\ldots,N\}}{\sigma\in\SG_N}}
P_{p,\sigma} [a,a_{\sigma(1)}\cdots a_{\sigma(p)}] b_1\cdots b_N
\otimes a_{\sigma(p+1)} \cdots a_{\sigma(N)} \otimes b
\end{equation}

$\bullet$ from $[\cR^{\prime 1,2}_N,r^{2,3}]$, we get
\begin{equation} \label{contr:2}
\sum_{\stackrel{p\in\{0,\ldots,N\}}{\sigma\in\SG_N}}
P_{p,\sigma} a_{\sigma(1)}\cdots a_{\sigma(p)} b_1\cdots b_N
\otimes [a,a_{\sigma(p+1)} \cdots a_{\sigma(N)}] \otimes b
\end{equation}

$\bullet$ from $[\cR^{\prime 1,3}_N,r^{2,3}]$, we get
\begin{equation} \label{contr:3}
\sum_{\stackrel{p\in\{0,\ldots,N\}}{\sigma\in\SG_N}}
P_{p,\sigma} a_{\sigma(1)}\cdots a_{\sigma(p)} b_1\cdots b_N
\otimes [a_{\sigma(p+1)},[ \cdots,[a_{\sigma(N)},a]] \otimes b
\end{equation}

$\bullet$ the only contribution from $[\cR^{\prime 2,3}_N,r^{1,3}]$ is  
the projection of $(\ref{contr:3})^{2,1,3}$ on 
$(U(\G)\otimes U(\A)\otimes U(\B)_0)^{*,*,1}_\univ$, which is zero since 
$N\neq 0$; all other contributions are zero. So we get 
$(\ref{contr:1})+(\ref{contr:2})+(\ref{contr:3}) = 0$. Selecting 
the terms where $a$ is in the first factor, we get   
$\sum_{\stackrel{p\in\{0,\ldots,N\}}{\sigma\in\SG_N}}
P_{p,\sigma} [a,a_{\sigma(1)}\cdots a_{\sigma(p)}] 
\otimes a_{\sigma(p+1)} \cdots a_{\sigma(N)}  = 0$
(equality in $\cF\cA_{N+1}^{\otimes 2}$), so $P_{p,\sigma} = 0$
unless $p=0$. Therefore $\cR'_N \in (U(\B)\otimes U(\A))_\univ$. 

We then get 
\begin{equation} \label{pre:Lie}
\sum_{\sigma\in\SG_N} P_{0,\sigma} [a_{\sigma(1)},\cdots,[a_{\sigma(N)},a]]
= \sum_{\sigma\in\SG_N} P_{0,\sigma} [a_{\sigma(1)}\cdots`a_{\sigma(N)},a]
\end{equation}
(equality in $\cF\cA_{N+1}$). Let us now show that this implies 
$P:= \sum_{\sigma\in\SG_N} P_{0,\sigma} a_{\sigma(1)}\cdots a_{\sigma(N)}\in\cF\cL_N$. 

Recall that $\cF\cA_{N+1}$ is the enveloping algebra $U(\cF\cL_{N+1})$, and
has therefore a Hopf lagebra structure. Denote by $\Delta_{\cF\cA_{N+1}}$ its
coproduct. The left side of (\ref{pre:Lie}) is primitive for thsi coproduct, 
so its right side is also primitive. We get 
$$
[\Delta_{\cF\cA_{N+1}}(P) - P \otimes 1 - 1 \otimes P,x\otimes 1 
+ 1 \otimes x] = 0. 
$$ 
Moreover, we have the identity $m([X,x\otimes 1 + 1\otimes X]) = 
[m(X),x]$, for any $X\in (\cF\cA_{N+1})^{\otimes 2}$, where 
$m : (\cF\cA_{N+1})^{\otimes 2} \to \cF\cA_{N+1}$ is the multiplication  
map. Therefore, 
$$
[m\big( \Delta_{\cF\cA_{N+1}}(P) - P \otimes 1 - 1 \otimes P\big) , x] = 0,  
$$
which implies that $m\big( \Delta_{\cF\cA_{N+1}}(P) - P \otimes 1 - 1 \otimes P\big)$
is a polynomial in $x$. Since this is also a homogeneous polynomial of
degree $N$ in $x_1,\ldots,x_N$, we get 
$$
m\big( \Delta_{\cF\cA_{N+1}}(P) - P \otimes 1 - 1 \otimes P\big) = 0. 
$$
Now is $\xx$ is any Lie algebra, the kernel of $U(\xx) \to U(\xx)$, 
$x\mapsto m(\Delta(x) - x \otimes 1   - 1 \otimes x)$ is equal to 
$\xx$; indeed, this map identifies via the symmetrisation 
isomorphism $U(\xx) \stackrel{\sim}{\to} S^\bullet(\xx)$
with $\sum_{i\geq 0} (2^i - 2) p_i$, where $p_i$ is the projector on 
$S^i(\xx)$. Therefore $P\in \cF\cL_{N+1}$, and since $P$ has degree 
$0$ in $x$ and $1$ in each generator $x_1,\ldots,x_N$, we get 
$P\in (\cF\cL_N)_{\sum_{i=1}^N \delta_i}$. Therefore 
$$
\cR'_N \in (U(\B)\otimes \A)_\univ. 
$$

Let us now project the identity 
$$
[\![r,\cR_N]\!] = \sum_{\stackrel{p,p',p''>0}{p+p'+p'' = N+1}} 
\cR^{2,3}_{p''}\cR^{1,3}_{p'}\cR^{1,2}_p - \cR^{1,2}_p\cR^{1,3}_{p'}
\cR^{2,3}_{p''}
$$
on $(U(\B)\otimes U(\A)^{\otimes 2})_\univ$. The projection of the 
right side is zero, since the third components of $\cR^{1,3}_{p'}$ and 
$\cR^{2,3}_{p''}$ are zero. On the other hand, the projection of 
$[\![r,\cR''_N]\!]$ on this space is zero for the same reason. 
Therefore the projection of  $[\![r,\cR'_N]\!]$ on 
$(U(\B)\otimes U(\A)^{\otimes 2})_\univ$ is zero. 
Since $[\cR^{\prime 1,2}_N,r^{1,3} + r^{2,3}]\in
(U(\G)\otimes U(\A) \otimes \B)_\univ$, 
$[\cR^{\prime 1,3}_N,r^{1,2}]\in
(U(\G)\otimes \B\otimes U(\G))_\univ$, and 
$[\cR^{\prime 2,3}_N,r^{1,2} + r^{1,3}]\in
(\A\otimes U(\G)^{\otimes 2})_\univ$, the projection of these
terms is also zero. So the projection of $[\![r,\cR_N]\!]$
coincides with that of $[\cR_N^{\prime 1,3},r^{2,3}]$. Then the 
projection of $\cR_N^{\prime 1,3}r^{2,3}$ is zero, and that of 
$$
r^{2,3}\cR_N^{\prime 1,3} = \sum_{\sigma\in\SG_N}
P_{0,\sigma} y_1\cdots y_N \otimes x \otimes 
y x_{\sigma(1)} \cdots x_{\sigma(N)}
$$    
is equal to $\sum_{\sigma\in\SG_N}
P_{0,\sigma} y_1\cdots y_N \otimes [x_{\sigma(1)},\ldots, 
[x_{\sigma(N)},x]] \otimes y$, which is therefore zero. 
Selecting the monomials containing $x$ as their last letter, 
we get  
$\sum_{\sigma\in\SG_N}
P_{0,\sigma} y_1\cdots y_N \otimes x_{\sigma(1)}\cdots 
x_{\sigma(N)} = 0$, so $\cR'_N = 0$. This proves Theorem 
\ref{thm:form:R}. \hfill \qed\medskip 

\subsection{A flatness criterion} \label{sect:flatness}

Let $(\A,\mu,\delta)$ be a finite dimensional Lie bialgebra
and let $\G = \A \oplus\A^*$ be its double Lie bialgebra. The inclusions 
$\A\hookrightarrow \G$ and $(\A^*)^{\cop} \hookrightarrow \G$ are Lie bialgebra
morphisms, where cop denotes the Lie bialgebra with the opposite coproduct. 
Then $\G$ is quasitriangular, 
with $r$-matrix $r_\A$. Assume that $(U(\G)[[\hbar]],m_0,\Delta,\cR_\A)$
is a quantization of $(\G,r)$, where the product of $U(\G)$ is not
deformed. Recall that this means that  $(U(\G)[[\hbar]],m_0,\Delta,\cR_\A)$
is a quasitriangular QUE algebra, such that 
$$
\big( (\cR_\A - 1)/\hbar \on{\ modulo\ } 
\hbar\big) = r_\A\in \G^{\otimes 2}. 
$$
Then $(U(\G)[[\hbar]],m_0,\Delta,\cR_\A)$ gives rise to the 
following quantized formal series Hopf algebras (QFSHA): 

$\bullet$ $(U(\G)[[\hbar]]^*,{}^t \Delta,{}^t m_0)$ is a quantization of
the FSHA $\cO(G) := U(\G)^*$; we denote it by $\cO_\hbar(G)$;   

$\bullet$ $U(\G)[[\hbar]]' \subset U(\G)[[\hbar]]$ is defined by 
$\hbar$-adic valuation conditions (see \cite{Dr:QG,Gav}) and is a 
quantization of the FSHA $\cO(G^*) := U(\G^*)^*$; we denote it by 
$\cO_\hbar(G^*)$. 

Before we state Theorem \ref{thm:flatness}, we recall the 
notion of a flat deformation of a morphism of vector spaces. 
If $f_0 : X_0 \to Y_0$ is a morphism of $\KK$-vector spaces, and 
$f : X \to Y$ is a morphism of topologically free 
$\KK[[\hbar]]$-modules, whose reduction is $f_0$, then 
$\Ker(f)$ is a divisible submodule of $X$, so 
$\Ker(f) / \hbar \Ker(f)$ injects in $X / \hbar X = X_0$, 
its image in $X_0$ being contained in $\Ker(f_0)$. 
 
\begin{lemma}  The following conditions are equivalent: 

1) $\Ker(f) / \hbar \Ker(f) = \Ker(f_0)$  

2) $\Imm(f)$ is divisible in $Y$.  

When these conditions are satisfied, we say that $f : X \to Y$
is a flat deformation of $f_0 : X_0 \to Y_0$. 
\end{lemma}

{\em Proof of Lemma.} Assume that $\Ker(f)/\hbar\Ker(f) = \Ker(f_0)$, 
then if $Z_0$ is a supplementary of $\Ker(f_0)$ in $X_0$, 
we have $X = \Ker(f)\oplus Z_0[[\hbar]]$. Then 
$\Imm(f)$ is the submodule of $Y$ generated by $f(Z_0)$. 
$f_0 : Z_0 \to Y_0$ is injective, so 
$f(Z_0) \cap \hbar Y = \{0\}$. Therefore $\Imm(f)$ is
divisible in $Y$. 

Conversely, if $\Imm(f)$ is divisible in $Y$, let us show that 
mod $\hbar : \Ker(f) \to \Ker(f_0)$ is surjective. Let 
$x_0\in \Ker(f_0)$ and let $x\in X$ be any lift of $x_0$. 
Then $f(x)\in \hbar Y$. So $f(x)\in \hbar Y \cap \Imm(f)$. 
Since $\Imm(f)$ is divisible, we have $f(x)\in \hbar \Imm(f)$, so 
$f(x) = \hbar f(x_1)$, where $x_1\in X$. So $f(x - \hbar x_1) = 0$. 
So $x-\hbar x_1\in \Ker(f)$, and its reduction modulo $\hbar$ is 
$x_0$. So mod $\hbar : \Ker(f) \to \Ker(f_0)$ is surjective. 
\hfill \qed\medskip 

Let us come back to the setup of the beginning of the section. 
If $\cA$ is a (Hopf) algebra, $\cA^\op$ denotes the algebra with 
the opposite product. 

\begin{thm} \label{thm:flatness}
1) The linear map $\ell  : U(\G)[[\hbar]]^{*,\op} \to U(\G)[[\hbar]]$, 
$\xi\mapsto \langle \cR_\A, \id\otimes\xi\rangle$ is a Hopf algebra 
morphism. Its image is contained in $\cO_\hbar(G^*)$. 

2) The reduction modulo $\hbar$ of $\ell : \cO_\hbar(G)^\op \to 
\cO_\hbar(G^*)$ is the FHSA morphism $\cO(G) \to \cO(G^*)$ dual to 
the composed map $G^* \twoheadrightarrow A^* \hookrightarrow G$
corresponding to the sequence of Lie bialgebra morphisms
$\G^* \twoheadrightarrow \A^* \hookrightarrow \G^\cop$. 

3) If we assume that $\cR_\A - 1 \in \hbar (U(\G) \otimes U(\G)\B)[[\hbar]]$,  
then $\ell : \cO_\hbar(G)^\op \to \cO_\hbar(G^*)$ is a flat deformation of 
$\cO(G) \to \cO(G^*)$, and $\Imm(\ell)$ is a quantization of $\cO(A^*)$.  
\end{thm}
 
{\em Proof.} Let us set $\delta_n = (\id - \eta\circ\eps)^{\otimes n}
\circ \Delta^{(n)}$. Then $\delta_n$ is a linear map 
$U(\G)[[\hbar]] \to U(\G)^{\otimes n}[[\hbar]]$, and 
$U(\G)[[\hbar]]'$ is defined as 
$$
\{x\in U(\G)[[\hbar]] \ |\  \forall n\geq 0, \ 
\delta_n(x)\in \hbar^n U(\G)^{\otimes n}[[\hbar]] \}. 
$$
Then the quasitriangular relations imply that 
for $\xi\in U(\G)[[\hbar]]^*$, 
$$
\delta_n(\ell(\xi))  = \langle (\cR_\A^{1,n+1}-1) 
\cdots (\cR_\A^{n,n+1}-1), \id^{\otimes n} \otimes\xi \rangle.  
$$
Since $\cR_\A - 1\in \hbar U(\G)^{\otimes 2}[[\hbar]]$, we get 
$\delta_n(\ell(\xi))\in \hbar^n U(\G)^{\otimes n}[[\hbar]]$, so 
$\ell(\xi)\in U(\G)[[\hbar]]'  = \cO_\hbar(G^*)$. This proves 1)
(see also \cite{KT}). 

Let us prove 2). $\ell$ mod $\hbar$ is a morphism of Poisson FSH 
algebras 
\begin{equation} \label{morphism}
\cO(G)^{\on{Poisson\ op}} \to \cO(G^*),
\end{equation} 
(``Poisson op'' means the Poisson 
algebra with opposite Poisson bracket), so it corresponds to a 
Lie bialgebra morphism $\G^*  \to \G^\cop$. This morphism is 
uniquely determined by the linear map $\mm_G / (\mm_G)^2 
\to \mm_{G^*} / (\mm_{G^*})^2$ induced by (\ref{morphism}),
using the identifications  $\mm_G / (\mm_G)^2 = \G^*$ and 
$\mm_{G^*} / (\mm_{G^*})^2 = \G$. 

For $\xi\in\G^*$, let $\wt\xi\in U(\G)^*$ be such that 
$\wt\xi(1) = 0$ and $(\wt\xi)_{|\G} = \xi$. Then 
$\wt\xi \in \mm_G$ and the class of $\wt\xi$ in  
$\mm_G / (\mm_G)^2 = \G^*$ is $\xi$. On the other hand, 
we have 
$\cR_\A \in 1 + \hbar r_\A + \hbar^2 U(\G)^{\otimes 2}[[\hbar]]$, 
so $\ell(\wt\xi) \in \hbar \langle r_\A,\id\otimes\xi \rangle 
+ \hbar^2 U(\G)[[\hbar]]$. Set $p_\A(\xi) = 
\langle r_\A,\id\otimes\xi \rangle$. Then $p_\A$ is the 
composed map $\G^* \twoheadrightarrow \A \hookrightarrow \G$.

So we have $\eps(\ell(\wt\xi)) = 0$, $\ell(\wt\xi)\in 
\cO_\hbar(G^*)$, and $\ell(\wt\xi) \in \hbar p_\A(\xi) + 
\hbar^2 U(\G)[[\hbar]]$. 
According to \cite{EGH}, Appendix, all this implies that 
$\big( \ell(\wt\xi) \on{\ mod\ }\hbar\big) \in \cO(G^*)$
belongs to $\mm_{G^*}$, and that its class in 
$\mm_{G^*} / (\mm_{G^*})^2 = \G$ is $p_\A(\xi)$. 
So the map $\mm_G / (\mm_G)^2 
\to \mm_{G^*} / (\mm_{G^*})^2$ induced by (\ref{morphism})
is $p_\A$. The morphism of Lie bialgebras corresponding to 
$\ell \on{\ mod\ }\hbar : \cO(G)\to\cO(G^*)$ is then $(p_\A)^t$, 
so it identifies with the composed map 
$\G^* \twoheadrightarrow \A^* \hookrightarrow \G^\cop$. 

Let us prove 3). Inclusion followed by multiplication induces 
a linear isomorphism 
\begin{equation} \label{isom:1}
U(\A)\otimes U(\A^*) \to U(\G).
\end{equation} 
Set $\cO(A) = U(\A)^*$, $\cO(A^*) = U(\A^*)^*$, then the dual map 
to (\ref{isom:1}) is a linear isomorphism  
\begin{equation} \label{isom:2}
\cO(G) \to \cO(A)\bar\otimes \cO(A^*), 
\end{equation} 
where $\bar\otimes$ is the tensor product of formal series rings. 

Let $\ell_0$ be the morphism $\cO(G) \to \cO(G^*)$ induced by the 
morphism $\G^*\to \G$. Then under (\ref{isom:2}), $\ell_0$ identifies 
with $\mm_A\bar\otimes\cO(A^*)$, where $\mm_A$ is the maximal ideal of 
$\cO(A)$. 

We now prove that the image by mod $\hbar : \cO_\hbar(G)^\op \to 
\cO(G)$ of $\Ker(\ell)$ is $\Ker(\ell_0)$. As a topological 
$\KK[[\hbar]]$-module, $\cO_\hbar(G) = U(\G)[[\hbar]]^*$ identifies 
with $\cO(G)[[\hbar]] = \cO(A)\bar\otimes\cO(A^*)[[\hbar]]$. 

Moreover, if $f\in\cO(A)$, $g\in\cO(A^*)$, $t\in U(\A)$ and $u\in U(\A^*)$,    
then $\langle tu, f\otimes g \rangle = \langle t,f\rangle \langle u,g \rangle$.
So if $f\in \cO(A)$ and $u\in U(\A^*)_0$, then 
$\langle tu,f\otimes 1 \rangle = \langle t,f \rangle \langle u,1 \rangle = 0$. 

So, if $f\in \mm_A$, then $\langle f,1 \rangle = 0$.    
Therefore, 
if $f\in \mm_A$, then $\langle \cR_\A, \id\otimes f\rangle = 0$. 
So $\Ker(\ell)$ contains $\mm_A[[\hbar]] \subset \cO(G)[[\hbar]] = 
\cO_\hbar(G)$. 

Moreover, $\Ker(\ell)$ is an ideal of $\cO_\hbar(G)^\op$. So 
$\Ker(\ell)$ contains the image of the map 
$\la : \mm_A[[\hbar]] \otimes \cO(B)[[\hbar]]\to \cO_\hbar(G)$
composed of the inclusions and the product of $\cO_\hbar(G)$.
The reduction modulo $\hbar$ of $\Imm(\la)$ is $\Ker(\ell_0)$. 
This proves that $(\on{mod\ }\hbar)(\Ker(\ell)) \supset \Ker(\ell_0)$. 
Since we also have the inverse inclusion, we get   
$(\on{mod\ }\hbar)(\Ker(\ell)) = \Ker(\ell_0)$. This proves that $\ell$
is a flat deformation of $\ell_0$. 

Therefore $\Imm(\ell)$ is a flat deformation of $\Imm(\ell_0) = \cO(A)$. 
It is also a Hopf subalgebra of $\cO_\hbar(G)$, so it is a QFSH algebra. 
Since $\cO(A)$ is a Poisson sub-FSHA of $\cO(G)$, $\Imm(\ell)$ is a 
quantization of $\cO(A)$. 
\hfill \qed\medskip

\subsection{Quantization of finite-dimensional Lie bialgebras}
\label{sect:fd}

Fix $\Phi\in\Assoc$. Set $\cR_\Phi = J_\Phi^{2,1}e^{t/2}(J_\Phi)^{-1}$. 
Then $\cR_\Phi\in(\wh U_2)^\times$. 
Set  $x^{1,2,\ldots,\{i,i+1\},\ldots,n;\Phi} = 
J_\Phi^{i,i+1} x^{1,2,\ldots,\{i,i+1\},\ldots,n}
(J_\Phi^{i,i+1})^{-1}$, for $x\in \wh U_n$. This operation  
satisfies bialgebra identities. We have  
$$
\cR_\Phi^{12,3;\Phi}  = 
\cR_\Phi^{1,3}\cR_\Phi^{2,3}, \ 
\cR_\Phi^{1,23;\Phi} = \cR_\Phi^{1,3}\cR_\Phi^{1,2}, 
$$
and 
$\sigma_{i,i+1}(x^{1,2,\ldots,\{i,i+1\},\ldots,n;\Phi}) = 
\cR_\Phi x^{1,2,\ldots,\{i,i+1\},\ldots,n;\Phi}
\cR_\Phi^{-1}$ for any $x\in \wh U_n$, and $\sigma_{i,i+1}(x) = 
x^{1,2,\ldots,i-1,i+1,i,i+2,\ldots,n}$. These are 
universal versions of the quasitriangular QUE algebra identities. 

Let $(\A,\mu,\delta)$ be a finite-dimensional Lie bialgebra. Let 
$\G$ be the double of $\A$, $r_\A$ be the $r$-matrix of $\G$. Set 
$\cR_{\Phi,\A} = \cR_\Phi(\hbar r_\A)$ and $\Delta_{\Phi,\A} = \Ad(J_\Phi(\hbar r_\A)) 
\circ \Delta_0$, where $\Delta_0$ is the coproduct of $U(\G)[[\hbar]]$. Then 
$$
(U(\G)[[\hbar]],m_0,\Delta_{\Phi,\A},\cR_{\Phi,\A})
$$
is a quasitriangular QUE algebra, and it is a quantization of $\G$. 
Since $\cR_{\Phi,\A}$ satisfies the hypothesis of Theorem 
\ref{thm:flatness}, 3), 
$\Imm(\ell)$ is a quantization of $\cO(A)$. It follows that $\Imm(\ell)^\circ$
is a quantization of $U(\A)$. Here   $\Imm(\ell)^\circ$ is the subspace of 
$\Imm(\ell)^*$ of all linear maps $\phi : \Imm(\ell) \to \KK[[\hbar]]$, 
such that the $\hbar$-adic valuation of $\phi(\mm^n)$ tends of infinity 
with $n$, where $\mm = \mm_0 + \hbar \Imm(\ell)$ and $\mm_0$ is the kernel of 
the counit of $\Imm(\ell)$. 

\medskip 

In order to formulate the propic version of this construction, we 
reformulate it as follows. 

We have a sequence of linear maps 
$$
\wh S^\bullet(\A)[[\hbar]] \stackrel{i}{\to} U(\G)[[\hbar]]^*
\stackrel{\ell}{\to} U(\G)[[\hbar]] \to  
\wh S^\bullet(\A)[[\hbar]] , 
$$
defined by 
$$
i : \wh S^\bullet(\A)[[\hbar]] \stackrel{\sim}{\to}  
S^\bullet(\A^*)[[\hbar]]^* \to  
(S^\bullet(\A)\bar\otimes S^\bullet(\A^*))[[\hbar]]^* 
\stackrel{\sim}{\to} 
U(\G)[[\hbar]]^*,  
$$
where the last map is dual to 
$U(\G) \to U(\A)\otimes U(\A^*) \stackrel{\Symm_\A^{-1} \otimes \Symm_{\A^*}^{-1}
}{\to} S^\bullet(\A)\otimes S^\bullet(\A^*)$, 
$\ell$ is $\xi\mapsto \langle \cR_\A,\id\otimes\xi\rangle$, and the third map
is 
$$
U(\G)[[\hbar]] \to U(\A)\otimes U(\A^*)[[\hbar]] 
\stackrel{\id\otimes\eps}{\to} U(\A)[[\hbar]] \stackrel{\Symm_\A^{-1}}{\to}
S^\bullet(\A)[[\hbar]]. 
$$ 
The composed map is an isomorphism $\la : \wh S^\bullet(\A)[[\hbar]] \to 
\wh S^\bullet(\hbar\A)[[\hbar]]$. Then for $x,y\in \wh S^\bullet(\A)[[\hbar]]$, we set 
\begin{equation} \label{zeroth:formula}
x * y = (\la^{-1}\circ p) \big( (\ell\circ i)(x)(\ell\circ i)(y)\big),   
\end{equation} 
\begin{equation} \label{first:formula}
\Delta(x) = (\la^{-1}\circ p)^{\otimes 2} 
\big( \Delta_{\Phi,\A}((\ell\circ i)(x))\big). 
\end{equation}
Then $(\wh S^\bullet(\A)[[\hbar]],*,\Delta)$ is the QFS algebra $\cO_\hbar(A)$. 

\subsection{Construction of quantization functors} \label{sect:inf:dim}    

\subsubsection{Definition of quantization functors}

In \cite{DQ}, we introduce various props: 

$\bullet$ $\underline{\on{Bialg}}$ (resp., $\underline{\on{Bialg}}_{\on{comm}}$)
is the prop of (commutative) bialgebras;   

$\bullet$ $\underline{\on{LBA}}$ (resp., $\underline{\on{LCA}}$)
is the prop of Lie bialgebras (resp., Lie coalgebras), with generators 
$(\mu,\delta)$ (resp., $\delta$). Both $\ul{\on{LCA}}$ and 
$\ul{\on{LBA}}$ are graded props, where $\delta$ has degree $1$
(resp., and $\mu$ has degree $0$), and $\wh{\ul{\on{LCA}}}$ and 
$\wh{\ul{\on{LBA}}}$ are the corresponding completions. 

\medskip 

There is a natural prop morphism $\ul{\on{Bialg}}_{\on{comm}} \to 
\wh S^\bullet(\wh{\ul{\on{LCA}}})$, 
which is a universal version of the functor $\C \mapsto (\wh S^\bullet(\C),m_0,\Delta)$, 
taking the Lie coalgebra $\C$ to a FSH algebra; when $\C$ is finite-dimensional, 
$\wh S^\bullet(\C) = U(\C^*)^*$. 

We also introduce the props $\underline{\on{Alg}}$, $\underline{\on{Poisson}}$
and $\underline{\on{Comm}}$ of associative, Poisson and commutative algebras. 
These props are naturally attached to operads. $\underline{\on{Alg}}$ has a natural 
filtration, and $\on{gr}(\underline{\on{Alg}}) = \underline{\on{Poisson}}$. 
There is a natural prop morphism $\underline{\on{Comm}} \to \underline{\on{Poisson}}$.  

Let $\ul{\on{Coalg}}$ be the prop of coalgebras. We show that the natural 
morphism 
$$
\bigoplus_{N\geq 0} \big(
\ul{\on{Alg}}(n,N) \otimes \ul{\on{Coalg}}(N,m)
\big)_{\SG_N} \to
\ul{\on{Bialg}}(n,m) 
$$
is a linear isomorphism. We introduce a filtration on $\ul{\on{Bialg}}$
by defining $F^i(\ul{\on{Bialg}}(n,m))$ as the image of 
$$
\bigoplus_{N\geq 0} \big(
F^i( \ul{\on{Alg}}(n,N)) \otimes \ul{\on{Coalg}}(N,m)
\big)_{\SG_N} 
$$
by this map. Then $\ul{\on{Bialg}}$ is a filtered prop. We set 
$\ul{\on{Bialg}}_{\on{Poisson}} = \gr(\ul{\on{Bialg}})$.  
We have a natural morphism  
$$
\nu : \ul{\on{Bialg}}_{\on{Poisson}} \to 
\wh S^\bullet(\wh{\ul{\LBA}}),
$$ 
which is a universal version of the 
assignment $\A\mapsto (\wh S^\bullet(\A),m_0,P,\Delta)$, where $P$
is the Poisson structure on the FSH algebra $\wh S^\bullet(\A)$. 

We have a commutative diagram of props 
$$
\begin{array}{ccc}
\ul{\on{Bialg}}_{\on{Poisson}}& \to & \wh S^\bullet(\wh{\ul{\LBA}})
\\ \downarrow & & \downarrow \\
\ul{\on{Bialg}}_{\on{comm}}& \to & \wh S^\bullet(\wh{\ul{\on{LCA}}}) , 
\\
\end{array}
$$
where $\ul{\on{Bialg}}_{\on{comm}}$ is the prop of commutative 
bialgebras.  

The props $\wh{\ul{\LBA}}$ and $\wh S^\bullet(\wh{\ul{\LBA}})$ are 
completions of graded props, so they can be viewed as filtered props. 

A quantization functor is a filtered prop morphism $\ul{\on{Bialg}}
\to \wh S^\bullet(\wh{\ul{\on{LBA}}})$, whose associated graded is 
$\nu$. A quantization functor is therefore the same as the data 
of $m\in \wh S^\bullet(\wh{\ul{\LBA}})(2,1)$ and 
$\Delta\in \wh S^\bullet(\wh{\ul{\LBA}})(1,2)$, satisfying 
bialgebra and quasiclassical limit relations. 

\subsection{Construction of quantization functors}

There exists a unique $\ul\pi\in \wh S^\bullet(\ul{\on{LA}})(2,1)$, 
such that for any Lie algebra $\G$, $(S^\bullet(\G),\ul\pi_{\G})$
is isomorphic to the algebra $U(\G)$, with isomorphism given by the 
inverse of $\Symm : S^\bullet(\G)\to U(\G)$. $\ul\pi$ can be derived from the 
Campbell-Baker-Hausdorff series. 

On the other hand, there exists $\ul\beta\in \wh S^\bullet
(\ul{\on{LBA}})(2,2)$, such that when $\A$ is a finite-dimensional 
Lie bialgebra, the map 
$$
\nu : S^\bullet(\A^*) \otimes S^\bullet(\A) 
\stackrel{\Symm_{\A^*}\otimes \Symm_\A}{\longrightarrow} 
U(\A^*) \otimes U(\A) \to U(\G) \to U(\A) \otimes U(\A^*)
\stackrel{\Symm_\A^{-1} \otimes \Symm_{\A^*}^{-1}}{\to} 
S^\bullet(\A) \otimes S^\bullet(\A^*) 
$$
(the middle maps are induced by injections followed by 
multiplication in $U(\G)$) 
is such that 
$$
\langle \nu(x\otimes \xi), \xi'\otimes x'\rangle = 
\langle \xi\otimes\xi', \beta_{\A}(x\otimes x')\rangle, 
$$
for any $x,x'\in S^\bullet(\A)$ and $\xi,\xi'\in S^\bullet(\A^*)$. 

All this implies that  there exists $\ul{D\mu}$ in 
$\wh S^\bullet(3,3)$ (we index the ``in'' entries by $(1,\bar 2,3)$
and the ``out'' entries by $(\bar 1,2,4)$) of the form 
$$
\ul{D\mu} = \ul{\pi}^{1b,\bar 1} \circ \ul{\beta}^{a3,b2}
\circ \ul{\kappa}^{\bar 2,a4}, 
$$ 
($a$ and $b$ are inner indices) where we defined
$\ul{\pi}\in \wh S^\bullet(\ul{\on{LA}})(2,1)$
and $\ul{\kappa}\in \wh S^\bullet(\ul{\on{LCA}})(1,2)$, such that if   
we express the product map $U(\G) \otimes U(\G) \to U(\G)$ 
as a map $m_\G : S^\bullet(\A)\otimes S^\bullet(\A^*) \otimes 
S^\bullet(\A)\otimes S^\bullet(\A^*)
\to S^\bullet(\A)\otimes S^\bullet(\A^*)$, then 
$$
\langle m_\G(x_1\otimes\xi_1 \otimes x_2\otimes\xi_2), 
\xi\otimes x'\rangle = 
\langle \ul{D\mu}_\A(x_1\otimes x'\otimes x_2), 
, \xi'\otimes \xi_1 \otimes\xi_2\rangle.  
$$

There exists $\al'\in (\wh S^\bullet(\A)\wh\otimes
\wh S^\bullet(\B)^{\wh\otimes 2})_\univ$, such that for any 
finite-dimensional Lie bialgebra $\A$, the composed map 
$$
S^\bullet(\A) \stackrel{\al}{\to} U(\G) 
\stackrel{(\rho_\A)^{-1}}{\to} S^\bullet(\A) \otimes S^\bullet(\A^*)
$$  
is given by $x\mapsto \langle \al', \id\otimes\id\otimes x \rangle$. 
Let $\ul\al\in \wh S^\bullet(\wh{\ul\LBA})(2,1)$ be the element 
corresponding to $\al'$. Then $\ul\al_0 := 
\ul\al( - \otimes 1)$ is an  
element of $\wh S^\bullet(\wh{\ul{\LBA}})(1,1)$. 
$\ul\al_0$ corresponds to the composed map 
$$
S^\bullet(\A) \stackrel{\al}{\to} U(\G) 
\stackrel{(\rho_\A)^{-1}}{\to} S^\bullet(\A) \otimes S^\bullet(\A^*)
\stackrel{\id\otimes\eps}{\to} S^\bullet(\A). 
$$

In the same way
as Theorem \ref{thm:flatness}, we prove that  
$\ul\al( - \otimes 1)$ is
invertible (this relies on the fact that $\cR_\Phi - 1$
belongs to $(U(\G) \otimes U(\G)\B)_\univ$). 

We define $m_\Phi\in \wh S^\bullet(\wh{\ul{\LBA}})(2,1)$
as follows. Set $\beta_0 = \beta(-\otimes 1)$, then 
$\beta_0\in \wh S^\bullet(\wh{\ul{\LBA}})(1,2)$. We set 
$$
m_\Phi = (\ul\al_0)^{-1} \circ \ul\pi \circ 
(\ul\al\otimes \id) \circ 
(\id\otimes (\ul\beta\circ\ul\al_0)).  
$$
Taking into account the equality $\ul\kappa(1) = 1\otimes 1$, 
$m_\Phi$ corresponds to the diagram  

\setlength{\unitlength}{1mm}
\begin{picture}(50,27)(-50,48)

\put(-53,63){$x_1\in$}\put(-53,53){$x_2\in$}


\put(-45,52){\line(1,0){11}}
\put(-34,52){\line(0,1){4}}
\put(-34,56){\line(-1,0){11}}
\put(-45,56){\line(0,-1){4}}
\put(-44,53){$S^\bullet(\A)$}

\put(-45,62){\line(1,0){11}}
\put(-34,62){\line(0,1){4}}
\put(-34,66){\line(-1,0){11}}
\put(-45,66){\line(0,-1){4}}
\put(-44,63){$S^\bullet(\A)$}


\put(-20,50){\line(1,0){11}}
\put(-9,50){\line(0,1){4}}
\put(-9,54){\line(-1,0){11}}
\put(-20,54){\line(0,-1){4}}
\put(-19,51){$S^\bullet(\A^*)$}

\put(-20,55){\line(1,0){11}}
\put(-9,55){\line(0,1){4}}
\put(-9,59){\line(-1,0){11}}
\put(-20,59){\line(0,-1){4}}
\put(-19,56){$S^\bullet(\A)$}

\put(-20,60){\line(1,0){11}}
\put(-9,60){\line(0,1){4}}
\put(-9,64){\line(-1,0){11}}
\put(-20,64){\line(0,-1){4}}
\put(-19,61){$S^\bullet(\A^*)$}

\put(-20,65){\line(1,0){11}}
\put(-9,65){\line(0,1){4}}
\put(-9,69){\line(-1,0){11}}
\put(-20,69){\line(0,-1){4}}
\put(-19,66){$S^\bullet(\A)$}


\put(5,55){\line(1,0){11}}
\put(16,55){\line(0,1){4}}
\put(16,59){\line(-1,0){11}}
\put(5,59){\line(0,-1){4}}
\put(6,56){$S^\bullet(\A^*)$}

\put(5,60){\line(1,0){11}}
\put(16,60){\line(0,1){4}}
\put(16,64){\line(-1,0){11}}
\put(5,64){\line(0,-1){4}}
\put(6,61){$S^\bullet(\A)$}


\put(30,52){\line(1,0){11}}
\put(41,52){\line(0,1){4}}
\put(41,56){\line(-1,0){11}}
\put(30,56){\line(0,-1){4}}
\put(31,53){$S^\bullet(\A^*)$}

\put(30,62){\line(1,0){11}}
\put(41,62){\line(0,1){4}}
\put(41,66){\line(-1,0){11}}
\put(30,66){\line(0,-1){4}}
\put(31,63){$S^\bullet(\A)$}

\put(41,64){\vector(1,0){4}}
\put(52,64){\vector(1,0){3}}
\put(46,63){$\ul\al_0^{-1}$}
\put(43,53){$\ni 1$}


\put(55,62){\line(1,0){11}}
\put(66,62){\line(0,1){4}}
\put(66,66){\line(-1,0){11}}
\put(55,66){\line(0,-1){4}}
\put(56,63){$S^\bullet(\A)$}

\put(67,63){$\ni m_\Phi(x_1\otimes x_2)$}


\put(-20,62){\vector(-2,1){4}}
\put(-24,64){\vector(2,1){4}}
\put(-34,64){\vector(1,0){5}}
\put(-28,64){$\ul\al$}

\put(-20,52){\vector(-2,1){4}}
\put(-24,54){\vector(2,1){4}}
\put(-34,54){\vector(1,0){5}}
\put(-28,54){$\ul\al$}

\put(-3,59){$\ul\beta$}
\put(-9,57){\vector(4,1){4}}
\put(-5,61){\vector(-4,1){4}}
\put(1,61){\vector(4,1){4}}
\put(5,57){\vector(-4,1){4}}

\put(22,65){$\ul\pi$}
\put(22,52){$\ul\kappa$}

\put(-9,66){\vector(1,0){30}}
\put(21,53){\vector(-1,0){30}}
\put(16,62){\vector(2,1){5}}
\put(21,54){\vector(-2,1){5}}

\put(25,66){\vector(4,-1){5}}
\put(30,53){\vector(-1,0){5}}
\end{picture}

\noindent
which is the analogue of formula (\ref{zeroth:formula}). 

To define $\Delta_\Phi \in \wh S^\bullet(\wh{\ul{\LBA}})(1,2)$, 
we make the following observation. When $x\in S^\bullet(\A)$, 
$\ell(x) = \langle \cR_{\Phi,\A},\id\otimes x\rangle$. Then 
\begin{equation} \label{good:formula}
\Delta_{\Phi,\A}(\ell(x)) = \langle \cR_{\Phi,\A}^{1,3} 
\cR_{\Phi,\A}^{2,3}, \id\otimes\id\otimes x\rangle. 
\end{equation}
Recall that $\rho_\A^{\otimes 2}\circ \Delta_{\Phi,\A}\circ \ell$
is a map $S^\bullet(\A) \to S^\bullet(\A)^{\otimes 2}\otimes 
S^\bullet(\A^*)^{\otimes 2}$. Since $\cR_\Phi^{1,3}\cR_\Phi^{2,3}$
lies in $(U(\G)^{\wh\otimes 3})_\univ$, there exists 
$\ul\delta_\Phi\in \wh S^\bullet(\wh{\ul{\on{LBA}}})(3,2)$, 
such that 
\begin{align*}
& \forall x,x'_1,x'_2\in S^\bullet(\A), \  
\forall \xi'_1,\xi'_2\in S^\bullet(\A^*), 
\\  
& \langle \rho_\A^{\otimes 2} \circ\delta_{\Phi,\A}\circ\ell(x), 
\xi'_1\otimes x'_1 \otimes \xi'_2\otimes x'_2 \rangle =   
\langle (\ul\delta_\Phi)_\A(x\otimes x'_1\otimes x'_2), \xi'_1\otimes\xi'_2
\rangle. 
\end{align*}
Then we set 
$$
\Delta_\Phi = (\ul\al_0 \otimes \ul\al_0)^{-1} \circ \ul\delta_\Phi
(-\otimes -\otimes 1\otimes 1). 
$$
We have $\Delta_\Phi \in \wh S^\bullet(\wh{\ul{\LBA}})(1,2)$. 
$\Delta_\Phi$ corresponds to the diagram

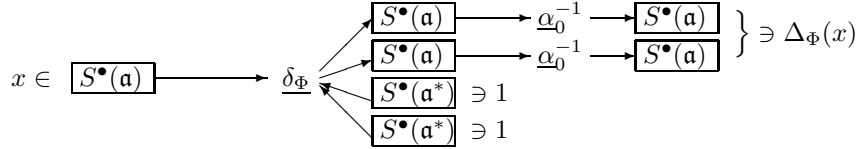
\setlength{\unitlength}{1mm}
\begin{picture}(50,27)(-50,48) 


\put(-53,58){$x\in$}


\put(-45,57){\line(1,0){11}}
\put(-34,57){\line(0,1){4}}
\put(-34,61){\line(-1,0){11}}
\put(-45,61){\line(0,-1){4}}
\put(-44,58){$S^\bullet(\A)$}

\put(-17,58){$\ul{\delta_\Phi}$}
\put(-34,59){\vector(1,0){15}}
\put(-12,60){\vector(1,1){7}}
\put(-12,59){\vector(3,1){7}}
\put(-5,56){\vector(-3,1){7}}
\put(-5,51){\vector(-1,1){7}}


\put(-5,50){\line(1,0){11}}
\put(6,50){\line(0,1){4}}
\put(6,54){\line(-1,0){11}}
\put(-5,54){\line(0,-1){4}}
\put(-4,51){$S^\bullet(\A^*)$}

\put(-5,55){\line(1,0){11}}
\put(6,55){\line(0,1){4}}
\put(6,59){\line(-1,0){11}}
\put(-5,59){\line(0,-1){4}}
\put(-4,56){$S^\bullet(\A^*)$}

\put(-5,60){\line(1,0){11}}
\put(6,60){\line(0,1){4}}
\put(6,64){\line(-1,0){11}}
\put(-5,64){\line(0,-1){4}}
\put(-4,61){$S^\bullet(\A)$}

\put(-5,65){\line(1,0){11}}
\put(6,65){\line(0,1){4}}
\put(6,69){\line(-1,0){11}}
\put(-5,69){\line(0,-1){4}}
\put(-4,66){$S^\bullet(\A)$}

\put(6,62){\vector(1,0){10}}
\put(17,61){$\ul{\al}_0^{-1}$}
\put(24,62){\vector(1,0){6}}

\put(6,67){\vector(1,0){10}}
\put(17,66){$\ul{\al}_0^{-1}$}
\put(24,67){\vector(1,0){6}}


\put(30,60){\line(1,0){11}}
\put(41,60){\line(0,1){4}}
\put(41,64){\line(-1,0){11}}
\put(30,64){\line(0,-1){4}}
\put(31,61){$S^\bullet(\A)$}

\put(30,65){\line(1,0){11}}
\put(41,65){\line(0,1){4}}
\put(41,69){\line(-1,0){11}}
\put(30,69){\line(0,-1){4}}
\put(31,66){$S^\bullet(\A)$}

\put(43,64){$\Big{\rbrace}\ni \Delta_\Phi(x)$}
\put(8,56){$\ni 1$}
\put(8,51){$\ni 1$}

\end{picture}

\noindent

One checks that $(m_\Phi,\Delta_\Phi)$ satisfy the bialgebra and 
quasiclassical limit relations, so $(m_\Phi,\Delta_\Phi)$ defines
a quantization functor of Lie bialgebras. 

\begin{remark}
If we use the propic analogue of formula (\ref{first:formula}) 
instead of 
(\ref{good:formula}), one would {\it a priori} get 
a diagram containing 
cycles (and therefore not suited for infinite-dimensional 
Lie bialgebras). In fact, this diagram concides with that of 
$\Delta_\Phi$ and does not contain cycles. 

More precisely, the {\it cyclization} of the prop $\ul{\LBA}$
is the prop $\ul{\LBA}_{\on{fd}}$, where we allow
cycles in the diagrams and the relations between generators
are the same as those of $\ul{\LBA}$. We have a prop morphism 
$\ul{\LBA} \to {\ul{\LBA}}_{\on{fd}}$ and therefore a linear map 
$$
\wh S^\bullet(\wh{\ul{\LBA}})(1,2) \to \wh S^\bullet(\wh{\ul{\LBA}}_{\on{fd}})(1,2).  
$$
The analogue of formula (\ref{first:formula}) 
yields an element $\wt\Delta_\Phi$ 
in this right side of this map, 
$\Delta_\Phi$ belongs to the left side of this map, and the image of 
$\Delta_\Phi$ is $\wt\Delta_\Phi$. 
\end{remark}

\newpage 

\begin{appendix}

\section{Commutation statements}

Recall that $m(r)$ be the element of $(\cF\cA_1 \otimes \cF\cA_1)_{\SG_1} 
\subset U(\G)_\univ$ equal to $x_1 \otimes y_1$. $m(r)$ is the universal 
version of $\sum_i a_i b_i$, where $r = \sum_i a_i \otimes b_i$. 
Moreover, if $(\G,\mu_\G,\delta_\G)$ is any double Lie algebra, the inner derivation 
$\ad(m(r_\A))$ coincides with the derivation $\mu_\G \circ \delta_\G$. 

\begin{prop} \label{m(r):central}
For any $Y \in (U(\G)^{\otimes n})_\univ$, we have 
\bn \label{dna}
[m(r)^1 + \cdots + m(r)^n, Y] = 0. 
\end{equation}
\end{prop}

{\em Proof.} Let us fix $N\geq 0$. We will prove (\ref{dna}) for any
$Y\in (U(\G)^{\otimes n})_{\univ,N}$. We have: 

(i) if $Y\in (U(\G)^{\otimes n})_{\univ,N}$ satisfies (\ref{dna}), 
and $1\leq i_1 < \cdots < i_n \leq m$, then so does 
$Y^{i_1,\ldots,i_n} \in (U(\G)^{\otimes m})_{\univ,N}$; 

(ii) if $Y\in (U(\G)^{\otimes n})_{\univ,N}$ satisfies (\ref{dna}), 
and $\sigma\in \SG_n$, then so does $Y^{\sigma(1),\ldots,\sigma(n)}$; 

(iii) let 
$$
m_{12} : (U(\G)^{\otimes n})_{\univ,N} \to (U(\G)^{\otimes n-1})_{\univ,N} 
$$   
be the map 
$$
((\cF\cA_N^{\otimes n})_{\sum_i \delta_i}
\otimes (\cF\cA_N^{\otimes n})_{\sum_i \delta_i})_{\SG_N}
\to 
((\cF\cA_N^{\otimes n-1})_{\sum_i \delta_i}
\otimes (\cF\cA_N^{\otimes n-1})_{\sum_i \delta_i})_{\SG_N}, 
$$
taking $\otimes_{i=1}^n (P_i \otimes Q_i)$ to 
$(P_1 P_2 \otimes Q_1 Q_2) \otimes \otimes_{i=3}^{n}(P_i \otimes Q_i)$. 
Then if $Y$ satisfies (\ref{dna}), so does $m_{12}(Y)$. 

(i) allows to reduce the proof of (\ref{dna}) to the case where 
$Y\in ((U(\G)_{>0})^{\otimes n})_{\univ,N}$ (the common intersection of all
partial counit maps $(\id^{\otimes i-1} \otimes \eps \otimes \id^{\otimes n-i-1})_\univ$).  
Then (ii) and (iii), together with the fact that $m_{12}$ is a surjective  
map 
$$
((U(\G)_{>0})^{\otimes n})_{\univ,N} \to ((U(\G)_{>0})^{\otimes n-1})_{\univ,N}, 
$$
allows to restrict the proof of (\ref{dna}) to the case where 
$Y\in (\G^{\otimes 2N})_{\univ,N}$. Now this space is spanned by all 
products $r^{i_1,j_1} \cdots r^{i_N,j_N}$, where $(i_1,\ldots,j_N)$
are such that $[1,2N] = \coprod_{\al = 1}^N \{i_\al,j_\al\}$ is a 
partition of $[1,2N]$. So it suffices to check (\ref{dna}) when 
$Y = r^{i_1,j_1} \cdots r^{i_N,j_N}$.  In this case, (\ref{dna}) follows from 
the facts that $\ad(m(r)^1 + \cdots + m(r)^n)$ is a derivation, and 
the explicit computation $[m(r)^1 + m(r)^2,r] = 0$. 
\hfill \qed \medskip 

\begin{prop} \label{t:invt}
If $X\in (U(\G)^{\otimes n})_\univ$, then 
$$
[t^{1,2},X^{12,3,\ldots,n+1}] = 0. 
$$
\end{prop}

{\em Proof.} We have 
$$
[t^{1,2},X^{12,3,\ldots,n+1}] = [m(r)^1,X]^{12,3,\ldots,n+1} 
- [m(r)^1 + m(r)^2, X^{12,3,\ldots,n+1}].  
$$
Applying Proposition \ref{m(r):central} to $X$ and taking the coproduct 
of the first factor, we get 
$$
[m(r)^1,X]^{12,3,\ldots,n+1} 
= - \sum_{k=3}^{n+1} [m(r)^k, X^{12,3,\ldots,n+1}], 
$$
and applying the same statement to $X^{12,3,\ldots,n+1}$, we get 
$$
[m(r)^1 + m(r)^2, X^{12,3,\ldots,n+1}]
= - \sum_{k=3}^{n+1} [m(r)^k, X^{12,3,\ldots,n+1}]. 
$$
Taking the difference of both equalities, we get Proposition 
\ref{t:invt}. 
\hfill \qed \medskip

\newpage 

\section{Computation of cohomology groups (proof of 
Proposition \ref{cohom:results})} \label{proof:coh:results}

Total antisymmetrization induces a linear isomorphism 
\bn \label{isom:complexes} 
\bigoplus_{k=1}^{n-1}
(\wedge^k(\A) \otimes \wedge^{n-k}(\B))_\univ \stackrel{\sim}{\to}
(\wedge^n(\G))_\univ. 
\end{equation}

Let us set 
$$
C^{p,q} = (\wedge^p(\A) \otimes \wedge^q(\B))_\univ 
$$
and define 
$$
\pa^{\prime (p,q)} : C^{p,q} \to C^{p+1,q}, 
$$
$$
x\mapsto \Alt_{1,\ldots,p+1}((\delta \otimes 
\id^{\otimes p+q-1})_\univ(x))
$$
and 
$$
\pa^{\prime\prime (p,q)} : C^{p,q} \to C^{p,q+1}, 
$$
$$
x\mapsto \Alt_{p+1,\ldots,p+q+1}((\id^{\otimes p} \otimes \delta \otimes 
\id^{\otimes q-1})_\univ(x))
$$
where $\Alt_{i_1,\ldots,i_k}$ is the partial antisymmetrization
in indices $(i_1,\ldots,i_k)$. We will sometimes write these maps simply
$\pa'$ and $\pa''$. 

Set $C^n = \oplus_{p,q\¦ p+q = n} C^{p,q}$, and define $\pa^n : 
C^n \to C^{n+1}$ as follows: the component $C^{p,q} \to C^{p',q'}$
of $\pa^n$ coincides with $\pa^{\prime (p,q)}$ if $(p',q') = (p+1,q)$, 
with $\pa^{\prime\prime (p,q)}$ if $(p'',q'') = (p,q+1)$ and is zero in 
all other cases. Then $(C^\bullet,\pa^\bullet)$ is a complex, and 
(\ref{isom:complexes}) sets up an isomorphism of $(C^\bullet,\pa^\bullet)$
with $( (\wedge^\bullet(\G))_\univ,\pa^\bullet)$. 

\subsection{Fine and total degrees}

If $\kk' = (k'_1,\ldots,k'_p)$ and $\kk'' = (k''_1,\ldots,k''_q)$
are collections of integers, such that 
$$
k'_1 + \cdots + k'_p = k''_1 + \cdots + k''_q = N,  
$$
let us denote by $(\A^{\otimes p} \otimes \B^{\otimes q})_{\univ,\kk',\kk''}$
the subspace of $(\A^{\otimes p} \otimes \B^{\otimes q})_{\univ}$
corresponding to the image of 
$$
\bigoplus_{ \stackrel{\eps'_1,\ldots,\eps'_p, 
\eps''_1,\ldots,\eps''_q\in \oplus_{i=1}^N \NN \delta_i \¦ }
{\¦\eps'_1\¦ = k'_1,\ldots,\¦\eps'_p\¦ = k'_p, 
\¦\eps''_1\¦ = k''_1,\ldots,\¦\eps''_q\¦ = k''_q} } 
\bigotimes_{i=1}^p (\cF\cL_N)_{\eps'_i}
\otimes \bigotimes_{j=1}^q (\cF\cL_N)_{\eps''_j}  
$$
by the projection map to the space of coinvariants. 
Then 
$$
(\A^{\otimes p} \otimes \B^{\otimes q})_\univ
= \bigoplus_{\kk',\kk''} (\A^{\otimes p} \otimes 
\B^{\otimes q})_{\univ,\kk',\kk''} .
$$

We say that the pair $(\kk',\kk'')$ is the "fine degree" of the component
$(\A^{\otimes p} \otimes \B^{\otimes q})_{\univ,\kk',\kk''}$. 
The "total degree" of this summand is $N = k'_1 + \cdots + k'_p = 
k''_1 + \cdots + k''_q$. If $Z\in (\A^{\otimes p} \otimes \B^{\otimes q})_\univ$, 
we denote by $Z_{\kk',\kk''}$ its component of fine degree $(\kk',\kk'')$. 

The differentials $\pa^n$ all have total degree $+1$. This implies that the 
cohomology groups $H^\bullet(C^\bullet,\pa^\bullet)$ are graded by the total 
degree. 

\subsection{Proof of $H^2(C^\bullet,\pa^\bullet) = 0$}

We have 
\begin{align*}
 H^2(C^\bullet,\pa^\bullet) = & \Ker\big(
\pa^{\prime(1,1)} : (\A\otimes \B)_\univ \to (\wedge^2(\A) \otimes \B)_\univ\big) 
\\ & 
\cap \Ker\big(
\pa^{\prime\prime(1,1)} : (\A\otimes \B)_\univ \to 
(\A \otimes \wedge^2(\B))_\univ\big) . 
\end{align*}
Moreover, $H^2(C^\bullet,\pa^\bullet) = \oplus_{N\geq 0} 
H^2(C^\bullet,\pa^\bullet)_N$, where each summand is the total degree
$N$ part of the cohomology group. Let us fix $N$ and compute this summand. 

Let $X\in (\A\otimes \B)_{\univ,N}$. In \cite{quant}, we proved that 
$$
(\pa^{\prime (1,1)}(X))_{(1,N),N+1} = [r^{1,3},X^{2,3}]. 
$$
Then $X$ can be recovered from $Y = (\pa^{\prime (1,1)}(X))_{(1,N),N+1}$
as follows: set 
$$
Y = \sum_{k=0}^N \sum_{\sigma\in \SG_{N}}
P_{k,\sigma} x \otimes x_1 \cdots x_N \otimes 
y_{\sigma(1)} \cdots y_{\sigma(k)} y 
y_{\sigma(k+1)} \cdots y_{\sigma(N)} ,  
$$ 
then 
$$
X = \sum_{\sigma\in \SG_N} P_{0,\sigma} x_1\cdots x_N \otimes 
y_{\sigma(1)} \cdots y_{\sigma(N)} . 
$$

In particular, $\pa^{\prime(1,1)}(X) = 0$ implies $X = 0$. 
Therefore $H^2(C^\bullet,\pa^\bullet)_N$, so $H^2(C^\bullet,\pa^\bullet) = 0$. 

\subsection{Proof of $H^3(C^\bullet,\pa^\bullet) = 
\KK \cdot [t^{1,2},t^{2,3}]$}

Direct computation shows that $H^3(C^\bullet,\pa^\bullet)_2 = 
\KK \cdot [t^{1,2},t^{2,3}]$.  

We will prove that for each integer $N>2$, the component 
total degree $N$ component $H^3(C^\bullet,\pa^\bullet)_N$ of 
$H^3(C^\bullet,\pa^\bullet)$ is zero. 

Let $(X^{aab},X^{abb}) \in (\wedge^2(\A) \otimes \B)_{\univ,N} 
\oplus (\A \otimes \wedge^2(\B))_{\univ,N}$ be such that 
$$
\pa' (X^{aab}) = 0 \quad \on{(equality\ in\ }
(\wedge^3(\A)\otimes \B)_\univ{)} 
$$
\bn \label{eq:second}
\pa'' (X^{aab}) + \pa' (X^{abb}) = 0 \quad \on{(equality\ in\ }
(\wedge^2(\A)\otimes \wedge^2(\B))_\univ{)} 
\end{equation}
$$
\pa'' (X^{abb}) = 0 \quad \on{(equality\ in\ }
(\A \otimes \wedge^3(\B))_\univ{)} . 
$$

If $Y\in (\A^{\otimes 2} \otimes \B^{\otimes 2})$, we denote 
by $Y_{*,*;(N,1)}$ the sum $\sum_{\kk'\¦ \¦\kk'\¦ = N+1} Y_{\kk',(N,1)}$. 
Then we have (see \cite{quant})
$$
(\pa''(X^{aab}))_{*,*;(N,1)} = [r^{1,4}+r^{2,4},(X^{aab})^{1,2,3}], 
$$ 
therefore 
$$
(\pa''(X^{aab}))_{(N,1);(N,1)} = [r^{1,4},(X^{aab})^{1,2,3}_{(N-1,1);N}].  
$$ 
In the same way, 
$$
(\pa'(X^{abb}))_{(N,1);(N,1)} = [r^{2,3},(X^{abb})^{1,3,4}_{N;(N-1,1)}],   
$$ 
so the $((N,1),(N,1))$ component of (\ref{eq:second}) gives 
$$
[r^{1,4},(X^{aab})^{1,2,3}_{(N-1,1);N}] + 
[r^{2,3},(X^{abb})^{1,3,4}_{N;(N-1,1)}] = 0,     
$$
which implies that for some $Y\in (\A\otimes \B)_{\univ,N}$, we have 
$$
(X^{aab})_{(N-1,1);N} = [r^{2,3},Y^{1,3}] 
\;\on{and}\;
(X^{abb})_{N;(N-1,1)} = -[r^{1,3},Y^{1,2}] . 
$$
Set $X' = X - \pa(Y)$.  Then $\pa(X') = 0$, and we have 
$$
(X^{\prime aab})_{(1,N-1);N} = (X^{\prime aab})_{(N-1,1);N} = 0, \; 
(X^{\prime abb})_{N;(1,N-1)} = (X^{\prime abb})_{N;(N-1,1)} = 0. 
$$
Then the only contribution to $(\pa'(X^{\prime aab}))_{*,*;(1,N)}$
is 
$$
(\pa'(X^{\prime aab}))_{*,*;(1,N)} = [r^{3,4},(X^{\prime aab})^{1,2,4}]. 
$$
This implies $X^{\prime aab} = 0$. In the same way, one proves 
$X^{\prime abb} = 0$. So $X' = 0$, and $X = \pa(Y)$. Therefore 
$H^3(C^\bullet,\pa^\bullet)_N = 0$, so $H^3(C^\bullet,\pa^\bullet) = 0$. 
\hfill \qed \medskip 

\end{appendix}


\frenchspacing

\begin{thebibliography}{DGK}

{\small

\bibitem{BN} D. Bar-Natan, On associators and the Grothendieck-Teichm\"uller 
group. I, Selecta Math. (N.S.) 4 (1998), 183-212. 


\bibitem{Dr:QG} V. Drinfeld, Quantum groups, Proceedings of ICM (Berkeley, 1986),
798-820, eds. AMS, Providence (RI), 1987. 

\bibitem{Dr:pb} V. Drinfeld, Some unsolved problems in the quantum groups 
theory, Lect. Notes Math. 1510 (1992), 1-8.


\bibitem{Dr:gal} V. Drinfeld, On quasitriangular quasi-Hopf algebras 
and a group closely connected with $\on{Gal}(\bar \QQ/\QQ)$, Leningrad Math. J. 
2:4 (1991), 829-60. 

\bibitem{quant} B. Enriquez, Quantization of Lie bialgebras
and shuffle algebras of Lie algebras, Selecta Math. (N.S.), 
7:3 (2001), 321-407. 

\bibitem{Enr:univ} B. Enriquez, One some universal algebras
associated to the category of Lie bialgebras, Adv. Math. 164:1
(2001), 1-23.  

\bibitem{DQ} B. Enriquez, P. Etingof, On the invertibility of 
quantization functors, in preparation.

\bibitem{EGH} B. Enriquez, F. Gavarini, G. Halbout, On the unicity of 
braidings of quasitriangular Lie bialgebras, math.QA/0207235. 

\bibitem{EK} P. Etingof, D. Kazhdan, Quantization of Lie 
bialgebras, I, Selecta Math. (N.S.) 2 (1996), 1-41.  

\bibitem{EK:2} P. Etingof, D. Kazhdan, Quantization of Lie 
bialgebras, II, Selecta Math. (N.S.) 4 (1998), 213-31. 

%

\bibitem{Liu} A. Liulevicius, Arrows, symmetries and 
representation rings, J. Pure. Appl. Algebra 19 (1980), 
259-73. 

\bibitem{Gav} F. Gavarini, Quantum duality principle, Ann. Inst. 
Fourier (Grenoble), 57:3 (2002), 809-34.  


\bibitem{KT} C. Kassel, V. Turaev, Biquantization of Lie bialgebras, 
Pacific J. Math. 195:2 (2000), 297-369. 




\bibitem{McLane} S. McLane, Categorical algebra, Bull. Amer. Math. Soc. 
71 (1965), 40-106.  


\bibitem{Pos} L. Positselski, letter to M. Finkelberg and R. Bezrukavnikov 
(in Russian), 1995. 

}
\end{thebibliography}
\end{document}